\documentclass[psamsfonts,12pt]{amsart}
\usepackage{mathpazo}
\usepackage{eulervm}

\usepackage[utf8]{inputenc}
\usepackage{amsmath}
\usepackage{mathtools}
\usepackage{amsthm}
\usepackage{amsfonts}
\usepackage{amssymb}
\usepackage{mathrsfs}
\usepackage[margin=1.0in]{geometry}
\usepackage[all]{xy}
\usepackage[english]{babel}
\usepackage{graphicx}
\usepackage{xparse,zref-savepos} 
\usepackage{tikz}
\usepackage{tikz-cd}
\usepackage{hyperref}
\usepackage[shortlabels]{enumitem}
\usetikzlibrary{matrix,arrows,decorations.pathmorphing,shapes.misc}
\usepackage{capt-of}
\usepackage[normalem]{ulem}
\usepackage{bbm}
\numberwithin{equation}{section}
\usepackage{courier}
\usepackage{enumitem}

\newtheorem{theorem}{Theorem}[section]
\newtheorem{lemma}[theorem]{Lemma}
\newtheorem{proposition}[theorem]{Proposition}

\newtheorem{corollary}[theorem]{Corollary}
\theoremstyle{definition}
\newtheorem{example}[theorem]{Example}

\newtheorem{definition}[theorem]{Definition}
\newtheorem{remark}[theorem]{Remark}

\newcommand{\eR}{\overline{\mathbb{R}}}

\newcommand{\R}{\mathbb{R}}
\newcommand{\N}{\mathbb{N}}

\newcommand{\met}{\mathbf{Met}}

\newcommand{\elp}[2]{\big\|#1\big\|_{#2}}
\newcommand\restr[2]{{
		\left.\kern-\nulldelimiterspace 
		#1 
		\vphantom{\big|} 
		\right|_{#2} 
}}

\newcommand{\cat}[1]{\mathbf{#1}}
\newcommand{\st}{\ \mid \ }
\newcommand{\isom}{\cong}
\newcommand{\xto}{\xrightarrow}
\newcommand{\incl}{\hookrightarrow}

\DeclareMathOperator{\Int}{Int}
\DeclareMathOperator{\Ind}{Ind}

\makeatletter
\def\moverlay{\mathpalette\mov@rlay}
\def\mov@rlay#1#2{\leavevmode\vtop{%
		\baselineskip\z@skip \lineskiplimit-\maxdimen
		\ialign{\hfil$\m@th#1##$\hfil\cr#2\crcr}}}
\newcommand{\charfusion}[3][\mathord]{
	#1{\ifx#1\mathop\vphantom{#2}\fi
		\mathpalette\mov@rlay{#2\cr#3}
	}
	\ifx#1\mathop\expandafter\displaylimits\fi}
\makeatother

\newcommand{\Expect}{{\rm I\kern-.3em E}}
\newcommand{\norm}[1]{\left\lVert#1\right\rVert}

\makeatletter
\@ifundefined{zsaveposx}{\let\zsaveposx\zsavepos}{}
\newcounter{hposcnt}
\renewcommand*{\thehposcnt}{hpos\number\value{hposcnt}}
\NewDocumentCommand{\lplabel}{o m}{%
	\stepcounter{hposcnt}%
	\zsaveposx{\thehposcnt l}%
	\zref@refused{\thehposcnt l}%
	\zref@refused{hpos0l}%
	\makebox[0pt][r]{\makebox[\dimexpr\zposx{\thehposcnt l}sp-\zposx{hpos0l}sp][l]{#2}}%
	\IfNoValueF{#1}
	{\def\@currentlabel{#2}\ltx@label{#1}}
}
\makeatother
\AtBeginDocument{\zsaveposx{hpos0l}}

\newlist{thmenum}{enumerate}{1}
\setlist[thmenum, 1]{label=(\arabic*), ref=\thetheorem ~(\arabic*)}
 
\begin{document}
	
\title{Universality of Persistence Diagrams and the Bottleneck and
  Wasserstein Distances}
\author{Peter Bubenik and Alex Elchesen}
\date{\today}
\maketitle

\begin{abstract}
We prove that persistence diagrams with the $p$-Wasserstein distance is the universal $p$-subadditive commutative monoid on an underlying metric space with a distinguished subset. This result applies to persistence diagrams, to barcodes, and to multiparameter persistence modules. In addition, the 1-Wasserstein distance satisfies Kantorovich-Rubinstein duality. 
\end{abstract}


\section{Introduction}

In computational settings persistent homology produces a persistence module
that is isomorphic to a finite direct sum of interval modules. The barcode and persistence diagram summarize this collection of intervals~\cite{Collins:2004,cohen2007stability}.
Distances between these summaries include the barcode metric based on the dissimilarity distance between intervals~\cite{Collins:2004, Carlsson:2005} and the bottleneck and Wasserstein distances for persistence diagrams based on the supremum-norm distance in the plane~\cite{cohen2007stability,cohen2010lipschitz}.
More generally, given a poset $P$ and a collection $\mathscr{S}$ of indecomposable persistence modules on $P$, one may consider persistence modules 
on $P$ that are isomorphic to finite direct sums of elements of $\mathscr{S}$. Examples include two-parameter persistence modules isomorphic to a finite direct sum of block modules~\cite{algstabzz,cochoy2016decomposition} and
multi-parameter persistence modules isomorphic to a finite direct sum of rectangle modules~\cite{bjerkevik2016stability}.

We study distances in a setting that includes all of these examples.
We start with a metric pair $(X,d,A)$ (Definition~\ref{def:metric-pair}). For persistence diagrams, the relevant metric pair is $(\R^2_{\leq},d,\Delta)$ or $(\overline{\R}^2_{\leq},d,\Delta)$, where $d$ is the metric obtained from one of the $q$-norms on $\R^2$ (Example~\ref{ex:half-plane}), and for barcodes it is $(\Int(\R),d,\emptyset)$, where $\Int(\R)$ denotes the set of intervals in $\R$ and $d$ is either the length of the symmetric difference or the Hausdorff distance (Example~\ref{ex:intervals}).
Given such a metric pair, we construct a free commutative monoid $(D(X,A),+,0)$ of \emph{persistence diagrams} on $(X,A)$ (Definition~\ref{def:persistence-diagram}) together with a family of \emph{Wasserstein distances} $W_p$ for all $p \in [1,\infty]$ (Definition~\ref{def:wasserstein}).
For persistence diagrams and barcodes we recover the metrics mentioned above (Example~\ref{ex:wasserstein}). 

We introduce the notion of \emph{$p$-subadditive commutative metric monoids}, an algebraic and metric object for discussing the above constructions. These are metric spaces which are also monoids and for which the metric is \emph{$p$-subadditive} (Definition \ref{def:metric-monoid}). We prove the following. For details see Definition \ref{def:overlineDp}.

\begin{proposition}[Functorial construction of persistence diagrams with Wasserstein distance] \label{prop:0}
  Given a metric pair $(X,d,A)$, $(D(X,A),W_p,+,0)$ is a $p$-subadditive commutative metric monoid and the canonical inclusion $i:(X,d,A) \to (D(X,A),W_p,0)$ is $1$-Lipschitz.
  Furthermore, this construction is functorial.
\end{proposition}

Our main result is the following. For details see Theorem \ref{thm:overlineDp}.

\begin{theorem}[Persistence diagrams with the Wasserstein distance as an adjoint] \label{thm:2}
The forgetful functor from $p$-subadditive commutative metric monoids to metric pairs has a left adjoint given by the functor in Proposition~\ref{prop:0}.
\end{theorem}

An equivalent statement of this result is that persistence diagrams with the Wasserstein distance are \emph{universal} as follows.



\begin{theorem}[Universality of persistence diagrams with the Wasserstein distance] \label{thm:2}
Given a metric pair $(X,d,A)$, $(D(X,A),W_p,+,0)$ is the universal $p$-subadditive commutative metric monoid obtained from $(X,d,A)$.
   That is, given any $p$-subadditive commutative metric monoid $(N,\rho,+,0)$ and $1$-Lipschitz map $\varphi:(X,d,A) \to (N,\rho,0)$, there is a unique $1$-Lipschitz monoid homomorphism $\tilde{\varphi}: (D(X,A), W_p,+,0) \to (N,\rho,+,0)$ such that
   $\tilde{\varphi} i = \varphi$.
 \begin{equation*}
    \begin{tikzcd}
      (X,d,A) \ar[r,"i"] \ar[dr,swap,"\varphi"] & (D(X,A),W_p,0)\\ 
      & (N,\rho,0)
    \end{tikzcd}
\ 
    \begin{tikzcd}
      (D(X,A),W_p,+,0) \ar[d,dashed,"\exists ^^21 \tilde{\varphi}"] \\
      (N,\rho,+,0)
    \end{tikzcd}
\ 
    \begin{tikzcd}
      (X,d,A) \ar[r,"i"] \ar[dr,swap,"\varphi"] & (D(X,A),W_p,0) \ar[d,dashed,"\tilde{\varphi}"]\\
      & (N,\rho,0)
    \end{tikzcd}
  \end{equation*}
\end{theorem}



From this it follows that the $p$-Wasserstein distance is the largest $p$-subadditive distance for persistence diagrams.
  For details see Theorem~\ref{thm:maximal} and Definition~\ref{def:dp}.

\begin{corollary}[Wasserstein distance as largest subadditive distance] \label{thm:1}
  Given a metric pair $(X,d,A)$, the $p$-Wasserstein distance
  $W_p$ is the largest $p$-subadditive metric on $D(X,A)$ compatible with $d$.
\end{corollary}


The following related result is of independent interest.
For details see Theorem~\ref{thm:CMon}.

\begin{theorem}[Symmetric monoidal structures for pointed metric spaces] \label{thm:4}
  For each $p \in [1,\infty]$ there is a symmetric monoidal category $\cat{Met_*^p}$ (Definitions \ref{def:p-strengthened} and \ref{def:xp} and   Corollary~\ref{cor:monoidal-category}) for which the category of commutative monoids internal to $\cat{Met_*^p}$ is the category of $p$-subadditive commutative metric monoids.
\end{theorem}

As a corollary to these results
we obtain Converse Stability Theorems (Section~\ref{sec:converse}).
When $p=1$, the Wasserstein distance satisfies Kantorovich-Rubinstein duality (Section~\ref{sec:duality}).
 
\subsection{Related work}

  Various metrics in applied topology have been shown to be universal in the sense that they are maximal stable distances.
  These include the matching distance~\cite{dAFL:10},
  the interleaving and bottleneck distance~\cite{lesnick2015theory}, the homology interleaving distance~\cite{blumberg2017universality}, and the Reeb graph edit distance~\cite{bauer2018reeb}.
  In contrast, we show that the Wasserstein distances are universal in the sense of category-theory (Theorem \ref{thm:2}).


    A version of our Corollary~\ref{thm:1} appears in~\cite{bubenik2018wasserstein}.
    Their version does not assume that the sum is finite but does restrict to the special case that the set $X$ is a set of objects in a Grothendieck category with local endomorphism rings and that the set $A$ consists of the zero object. Also, they do not show that their result follows from a functorial construction.


	The Wasserstein distances between persistence diagrams has been studied extensively by Divol and Lacombe \cite{divol2019understanding}. There they relate the Wasserstein distance between persistence diagrams to the classical Wasserstein distance on probability measures. Among other things, this allows for a version of Kantorovich-Rubinstein duality to be recovered for persistence diagrams. 

Note that all of the persistence diagrams defined in the present paper are finite by definition. Extensions of the ideas presented here to countable persistence diagrams, such as those in \cite{cohen2010lipschitz} \and \cite{Blumberg:2014}, and signed persistence diagrams, as well as to the setting of Radon measures, can be found in the sequel \cite{bubenik2020virtual}.

Skraba and Turner~\cite{skraba2020wasserstein} have shown that the Wasserstein distance for persistence diagrams of weighted cell complexes is stable.

\section{Background}\label{sec:backgroud}

\subsection{Metric spaces}

In order to include various distances arising in persistent homology,
we will use a less restrictive notion of metric space than is standard.
Such relaxed metrics are referred to as extended pseudometrics in the literature, but we will refer to them as metrics for brevity.

\begin{definition} \label{def:metric}
  A \emph{metric space} is a tuple $(X,d)$ where $X$ is a set and $d:X \times X \to [0,\infty]$ is a function satisfying
  $d(x,x) = 0$ for all $x \in X$ (point triviality), $d(x,y) = d(y,x)$ for all $x,y \in X$ (symmetry), and $d(x,y) \leq d(x,z) + d(z,y)$ for all $x,y,z \in X$ (triangle inequality).
  Given metric spaces $(X,d_X)$ and $(Y,d_Y)$, a \emph{metric map} is a function $f:X \to Y$ such that $d_Y(f(x),f(y)) \leq d_X(x,y)$ for all $x,y \in X$. 
  Given a set $X$, a metric space $(Y,d)$, and a function $g:X\to Y$, the \emph{pullback of $d$ along $g$}, denoted $g^*d$, is the metric on $X$ defined by $g^*d(x,x') := d(g(x),g(x'))$.
\end{definition}

\subsection{Monoids}
\label{sec:monoids}

A \emph{commutative monoid} $M = (M,+,0)$ is a set $M$ together with an associative commutative binary operation $+:M\times M\to M$ for which there exists an element $0\in M$ satisfying $m+0 = m$ for all $m\in M$, 
called the \emph{identity} element.
A \emph{monoid homomorphism} between commutative monoids $M = (M,+_M,0_M)$ and $N = (N,+_N,0_N)$ is a map $f:M\to N$ such that 
$f(a +_M b) = f(a) +_N f(b)$ for all $a,b\in M$ and
$f(0_M) = 0_N$.
A subset $P \subset M$ is a \emph{submonoid} if it contains $0$ and $+$ restricts to a binary operation on $P$.

Given a set $X$, the \emph{free commutative monoid on $X$}, denoted $D(X)$, is the set of all (finite) 
\emph{formal sums} of elements of $X$, with the monoid operations being given by addition of formal sums.
That is, $D(X)$ is the set of all functions $f:X\to \N\cup \{0\}$ with finite support and with the monoid operation given by the pointwise addition of functions.
Formal sums are also called {finite multisets}.
For $x \in X$, let $1_x:X \to \N \cup \{0\}$ be given by $1_x(x) = 1$ and $1_x(y)=0$ for all other $y \in X$.
As is customary, we denote $1_x$ by $x$.
With this convention, we may write any formal sum $\alpha \in D(X)$ as $\alpha = x_1 + \cdots + x_n$, where $n \geq 0$ and $x_1,\ldots,x_n \in X$.
We define the \emph{canonical inclusion} $i:X\to D(X)$ by $i(x) = x$.
%
  
An equivalence relation $\sim$ on a commutative monoid $M$ is called a \emph{congruence} if $a\sim b$ and $c\sim d$ implies $a+c \sim b + d$. If $\sim$ is a congruence then there is a well-defined commutative monoid structure on the set of equivalence classes $M/\!\!\sim$ given by $[a] + [b] := [a+b]$. 
%
Let $M$ be a commutative monoid and $P \subseteq M$ any submonoid. Define a relation $\sim$ on $M$ by 
\[a\sim b \iff \textup{ $\exists x,y\in P$ such that $a+x = b+y$}.\]
Then $\sim$ is a congruence and we denote the commutative monoid $M/\!\!\sim$ by $M/P$ and refer to it as the \emph{quotient of $M$ by $P$}. 

  A \emph{set pair} is a pair $(X,A)$ where $X$ is a set and $A$ is a nonempty subset of $X$.
  A map of pairs $f:(X,A) \to (Y,B)$ is a function $f:X \to Y$ such that $f(A) \subset B$.
  A \emph{pointed set} is a pair $(X,\{x_0\})$, which is denoted $(X,x_0)$.
  Given pointed sets $(X,x_0)$ and $(Y,y_0)$, a \emph{pointed function} $f:(X,x_0) \to (Y,y_0)$ is a function $f:X \to Y$ such that $f(x_0) = y_0$.

\begin{definition} \label{def:persistence-diagram}
  Given a pair $(X,A)$, let $D(X,A)$ denote the quotient monoid $D(X)/D(A)$.
  We call $D(X,A)$ the commutative monoid of \emph{persistence diagrams} on $(X,A)$.
Note that $D(X)/D(A) \isom D(X \setminus A)$.
Given a map of pairs $f:(X,A) \to (Y,B)$, there is an induced monoid homomorphism $f_*:D(X,A) \to D(Y,B)$ given by $f_*(x_1 + \cdots + x_n) = f(x_1) + \cdots + f(x_n)$.
Note that this also defines a pointed function $f_*:(D(X,A),0) \to (D(Y,B),0)$.
\end{definition}

\subsection{$p$-norms}

Let $p\in [1,\infty]$ and
$\mathbf{x} = (x_1,\dots,x_n)\in \R^n$.
For $p < \infty$, the
\emph{$p$-norm} \emph{of $\mathbf{x}$} is defined by
$\|\mathbf{x}\|_p = \left(\sum_{k = 1}^n |x_k|^p\right)^{1/p}$ and
for $p = \infty$, it is defined by
$\|\mathbf{x}\|_\infty = \max_{1\leq k\leq n}|x_k|$.
%
%
For $x = (x_1,\ldots,x_m)$, $y=(x_{m+1},\ldots,x_{m+n})$, and $z = (x_1,\ldots,x_{m+n})$, $\lVert(\norm{x}_p,\norm{y}_p)\rVert_p = \norm{z}_p$.
By the \emph{$\ell^p$-distance on $\R^n$} we mean the metric induced by the $p$-norm, i.e., 
$\|\mathbf{x} - \mathbf{y}\|_p$.
The fact that each $\|\cdot\|_p$ is a norm relies on the
Minkowski inequality:
for all $p\in [1,\infty]$ and $\mathbf{x},\mathbf{y}\in \R^n$,
$\|\mathbf{x} + \mathbf{y}\|_p \leq \|\mathbf{x}\|_p + \|\mathbf{y}\|_p$.
The $p$-norms are related as follows, which shows in particular
that the $p$-norms are decreasing in $p$.
%
For $\mathbf{x}\in \R^n$ and $1\leq p\leq q\leq \infty$ we have
$\|\mathbf{x}\|_q \leq \|\mathbf{x}\|_p\leq n^{\frac{1}{p}-
  \frac{1}{q}}\|\mathbf{x}\|_q$, and these inequalities are
attained. Here we adopt the convention $\frac{1}{\infty} = 0$.

Let $\eR$ denote the set of extended real numbers $[-\infty,+\infty]$.
  The $\ell^p$-distance on $\R^n$ extends to $\eR^n$, with the understanding that it may take a value of $\infty$.
  
 \subsection{Basic category theory}\label{sec:category_theory} A \emph{category} $\cat{C}$ consists of a class $\textup{obj}(\cat{C})$ of \emph{objects}, and for each pair of objects $X,Y\in \textup{obj}(\cat{C})$, a set $\cat{C}(X,Y)$ of \emph{morphisms} (or \emph{arrows}). The class of all morphisms of $\cat{C}$ is denoted $\textup{Hom}(\cat{C})$. A morphism $f\in \cat{C}(X,Y)$ is often denoted $f:X\to Y$.  We will often simply write $X\in\cat{C}$ to indicate that $X$ is an object of $\cat{C}$. A category is \emph{small} if $\textup{obj}(\cat{C})$ is a set as opposed to a proper class.
 
 The objects and morphisms of a category $\cat{C}$ are required to satisfy the following axioms. For any objects $X,Y,Z\in \cat{C}$ and morphisms $f\in \cat{C}(X,Y)$, $g\in\cat{C}(Y,Z)$, there exists a morphism $g\circ f\in \cat{C}(X,Z)$, called the \emph{composition of $f$ and $g$}. That is, there is a function $\cat{C}(X,Y)\times \cat{C}(Y,Z)\to \cat{C}(X,Z)$ given by $(f,g)\mapsto g\circ f$. We will often omit the $\circ$, writing $gf$ instead of $g\circ f$. Composition must be \emph{associative}, meaning that $(hg)f = h(gf)$ whenever this composition is defined. Finally, for all $X\in \cat{C}$, there exists a morphism $\textup{id}_X:X\to X$ such that, for all $W,Y\in\cat{C}$ and $f:W\to X$, $g:X\to Y$, we have $\textup{id}_Xf = f$ and $g\textup{id}_X = g$.
 
 A \emph{subcategory} $\cat{D}$ of $\cat{C}$ consists of a subclass $\textup{obj}(\cat{D})$ of $\textup{obj}(\cat{C})$ and a subclass $\textup{Hom}(\cat{D})$ of $\textup{Hom}(\cat{C})$ such that if $f:X\to Y\in \textup{Hom}(\cat{D})$ then $X,Y\in \textup{obj}(\cat{D})$, $\textup{id}_X\in \textup{Hom}(\cat{D})$ for all $X\in \textup{obj}(\cat{D})$, and $fg \in \textup{Hom}(\cat{D})$ whenever $f,g\in \textup{Hom}(\cat{D})$ and this composition is defined in $\cat{C}$. This definition guarantees that $\cat{D}$ is a category in its own right. $\cat{D}$ is a \emph{full subcategory} of $\cat{C}$ if $\cat{D}(X,Y) = \cat{C}(X,Y)$ for all $X,Y\in\cat{D}$.
 
 Objects $X,Y\in\cat{C}$ are said to be \emph{isomorphic} if there exists morphisms $f:X\to Y$ and $g:Y\to X$ such that $gf = \textup{id}_X$ and $fg = \textup{id}_Y$.
 \begin{example}
 	Let $\cat{Set}$ denote the category whose objects are sets and whose morphisms are functions between sets. Composition is given by the composition of functions and the identity morphism on a set $S$ is the identity function on $S$. Isomorphisms in $\cat{Set}$ are bijective functions.
 \end{example}
 \begin{example}
 	Let $\cat{Met}$ denote the category whose objects are metric spaces $X = (X,d)$ and whose morphisms are metric maps (see Definition \ref{def:metric}). Composition of metric maps is given by the composition of functions and the identity morphism on $X$ is the identity function on $X$. It is easily checked that the composition of metric maps is again a metric map, as is the identity function. Isomorphisms in $\cat{Met}$ are isometries.
 \end{example}

\begin{example}
	Let $\cat{CMon}$ denote the category whose objects are commutative monoids $M = (M,+,0)$ and whose morphisms are monoid homomorphisms. Composition of monoid homomorphisms is given by the composition of functions and the identity morphism on $M$ is the identity function on $M$. Isomorphisms in $\cat{CMon}$ are monoid isomorphisms.
\end{example}

A \emph{covariant functor} $F:\cat{C}\to\cat{D}$ consists of a map $F:\textup{obj}(\cat{C})\to \textup{obj}(\cat{D})$ and, for each $X,Y\in \cat{C}$, a map $F:\cat{C}(X,Y)\to \cat{D}(F(X),F(Y))$. For all $X,Y,Z\in \cat{C}$ and $f\in\cat{C}(X,Y)$, $g\in \cat{C}(Y,Z)$, these maps must satisfy $F(gf) = F(g) F(f)$ and $F(\textup{id}_X) = \textup{id}_{F(X)}$.

 A \emph{contravariant functor} $F:\cat{C}\to\cat{D}$ consists of a map $F:\textup{obj}(\cat{C})\to \textup{obj}(\cat{D})$ and, for each $X,Y\in \cat{C}$, a map $F:\cat{C}(X,Y)\to \cat{D}(F(Y),F(X))$. For all $X,Y,Z\in \cat{C}$ and $f\in\cat{C}(X,Y)$, $g\in \cat{C}(Y,Z)$, these maps must satisfy $F(gf) = F(f) F(g)$ and $F(\textup{id}_X) = \textup{id}_{F(X)}$. Note that a contravariant functor reverses the direction of arrows in the sense that if $f:X\to Y$ then $F(f):F(Y)\to F(X)$. 

\begin{example}
	Any category admits an identity functor $1_\cat{C}:\cat{C}\to\cat{C}$ which maps objects and morphisms to themselves.
\end{example}
\begin{example}
	If $\cat{D}$ is a subcategory of $\cat{C}$, then the inclusion $\cat{D}\hookrightarrow \cat{C}$ is a covariant functor.
\end{example}
\begin{example}
	Let $G:\cat{D}\to \cat{C}$ be a covariant functor. For each $C\in \cat{C}$, there is a covariant functor $\cat{C}(C,G-):\cat{D}\to \cat{Set}$ which sends $D\in \cat{D}$ to $\cat{C}(C,GD)\in \cat{Set}$, and which sends a morphism $f:D\to D'$ in $\cat{D}$ to the set map $\cat{C}(C,Gf):\cat{C}(C,GD)\to \cat{C}(C,GD')$ given by $g:C\to GD \mapsto (Gf)g:C\to GD'$.
	
	Similarly, given a covariant functor $F:\cat{C}\to\cat{D}$ and $D\in \cat{C}$, there is a contravariant functor $\cat{D}(F-,D):\cat{C}\to \cat{Set}$ which sends $C\in \cat{C}$ to $\cat{D}(FC,D)\in \cat{Set}$, and which sends a morphism $f:C\to C'$ in $\cat{C}$ to the set map $\cat{D}(Ff,D):\cat{D}(FC',D)\to\cat{D}(FC,D)$ given by $h:FC'\to D \mapsto h(Ff):FC\to D$.
	
	As a special case, if $G = 1_\cat{C}$, for any $C\in \cat{C}$ we obtain the \emph{covariant hom-functor} $\cat{C}(C,-)$. Similarly, if $F= 1_\cat{C}$ then for any $C\in\cat{C}$ we obtain the \emph{contravariant hom-functor} $\cat{C}(-,C)$.
\end{example}

\begin{example}
	The \emph{forgetful functor} $U:\cat{CMon}\to\cat{Set}$ sends a commutative monoid to its underlying set, and sends a monoid homomorphism to the function defining it. That is, $U$ sends a commutative monoid to the set obtained by ``forgetting" the monoid structure.
	
	Similarly, there is a forgetful functor $U:\cat{Met}\to\cat{Set}$.
\end{example}

\begin{example}
	The \emph{free commutative monoid functor} $D:\cat{Set}\to\cat{CMon}$ sends a set $X$ to the free commutative monoid $D(X)$ (see Section \ref{sec:monoids}) and sends a function $f:X\to Y$ the the monoid homomorphisms $D(f):D(X)\to D(Y)$ given by $x_1 + \dots + x_n \mapsto f(x_1) + \dots + f(x_n)$.
\end{example}

A \emph{natural transformation} $\alpha:F\Rightarrow G$ between functors $F,G:\cat{C}\to \cat{D}$ (either both covariant or both contravariant), denoted $\alpha:F\Rightarrow G$, consists of, for each $C\in \cat{C}$, a morphism $\alpha_C:F(C)\to G(C)$. If $F$ and $G$ are both covariant then we require that for any $D,E\in \cat{C}$ and any morphism $f\in \cat{C}(D,E)$, we have $\alpha_E F(f) = G(f)\alpha_D$. If $F$ and $G$ are both contravariant then we require that $\alpha_D F(f) = G(f) \alpha_E$. The morphism $\alpha_C$ is called \emph{the component of $\alpha$ at $C$}. If all of the components of $\alpha$ are isomorphisms, then $\alpha$ is called a \emph{natural isomorphism}. If $\alpha$ is a natural isomorphism from $F$ to $G$ then we say that $F$ and $G$ are \emph{naturally isomorphic} and write $F\cong G$.

An \emph{adjunction} between categories $\cat{C}$ and $\cat{D}$ consists of functors $F:\cat{C}\to \cat{D}$ and $G:\cat{D}\to\cat{C}$ such that, for each $D\in \cat{D}$ there is a natural isomorphism $\cat{D}(F-,D)\cong \cat{C}(-,GD)$ and for each $C\in \cat{C}$ there is a natural isomorphism $\cat{D}(FC,-)\cong \cat{C}(C,G-)$. In this case, we say that $F$ is the \emph{left adjoint of $G$} and \emph{$G$ is the right adjoint of $F$}. We also write $F\dashv G$.


\subsection{Universal properties and adjunctions}
\label{sec:univ_props_adjunctions}


Consider a functor $U: \cat{D} \to \cat{C}$ between categories $\cat{D}$ and $\cat{C}$ and let $X \in \cat{C}$.
An object $Y \in \cat{D}$ satisfies a \emph{universal property} with \emph{universal element} $i: X \to UY \in \cat{C}$ if for every object $Z \in \cat{D}$ and morphism $f: X \to UZ$ in $\cat{C}$ there is a unique morphism $g:Y \to Z$ in $\cat{D}$ such that $Ug \circ i = f$.
\begin{equation} \label{eq:universal-property}
  \begin{tikzcd}
    X \ar[r,"i"] \ar[dr,swap,"f"]& UY \\
    & UZ
  \end{tikzcd}
  \quad \quad \quad
  \begin{tikzcd}
    Y \ar[d,dashed,"\exists ^^21 g"] \\
    Z
  \end{tikzcd}
  \quad \quad \quad
  \begin{tikzcd}
    X \ar[r,"i"] \ar[dr,swap,"f"]& UY \ar[d,dashed,"Ug"]\\
    & UZ
  \end{tikzcd}
\end{equation}

For example, consider the forgetful functor $U:\cat{CMon} \to \cat{Set}$ and let $X \in \cat{Set}$.
Then the free commutative monoid $D(X) \in \cat{CMon}$ satisfies a universal property with universal element the canonical inclusion $i:X \to UD(X) \in \cat{Set}$.

Given $C\in\cat{C}$ and a functor $U:\cat{D}\to\cat{C}$, the \emph{comma category} $C\downarrow U$ is the category whose objects are pairs $(D,f:C\to UD)$, where $D\in\cat{D}$ and $f:C\to UD$ is a morphism in $\cat{C}$, and whose morphisms $(D,f:C\to UD) \to (D',f':C\to UD')$ are morphisms $g:D\to D'$ in $\cat{D}$ for which $Ug\circ f = f'$. Then a more succinct way of stating \eqref{eq:universal-property} is that $(Y,i)$ is the initial object of the category $X \downarrow U$.

\begin{remark}The \emph{category of elements} of a functor $F:\cat{C}\to\cat{Set}$ is the category whose objects are pairs $(C,x)$, where $C\in\cat{C}$ and $x\in FC$, and whose morphisms $(C,x)\to (C',x')$ are morphisms $f:C\to C'$ in $\cat{C}$ for which $Ff(x) = x'$. Then $X \downarrow U$ is precisely the category of elements of the functor $\cat{C}(X,U-): \cat{D} \to \cat{Set}$. Then yet another way of stating \eqref{eq:universal-property} is that there is a natural isomorphism $\cat{C}(X,U-) \isom \cat{D}(Y,-)$ defined by $i \in \cat{C}(X,UY)$.
That is, $\cat{C}(X,U-)$ is represented by $Y$ via $i$.
\end{remark}

The following lemma shows how a family of universal properties can be used to obtain an adjunction. We will use it frequently throughout.
\begin{lemma}[{\cite[Lemma 4.6.1]{riehl2017category}}]\label{lem:adjoint}
  A functor $U: \cat{D} \to \cat{C}$ admits a left adjoint
  if and only if for each object $X$ in $\cat{C}$ the comma category $X \downarrow U$ has an initial object.
\end{lemma}

For example, the forgetful functor $U:\cat{CMon} \to \cat{Set}$ has a left adjoint $D: \cat{Set} \to \cat{CMon}$, the free commutative monoid functor. This follows from Lemma \ref{lem:adjoint} since for any set $X$, $X\downarrow U$ has the initial object $(D(X),i_X)$, where $D(X)$ is the free commutative monoid and $i_X:X\to D(X)$ is the canonical inclusion (see Section \ref{sec:monoids}).


  A special case of \eqref{eq:universal-property} has $U: \cat{D} \to \cat{C}$ being the inclusion of a full subcategory.
\begin{definition}
  A \emph{reflective} subcategory of a category $\cat{C}$ is a full subcategory $\cat{D}$ of $\cat{C}$ such that the inclusion $\cat{D} \hookrightarrow \cat{C}$ has a left adjoint. The left adjoint of the inclusion is called the \emph{reflector}.
\end{definition}

  \subsection{Symmetric monoidal categories and internal objects}\label{sec:sym_mon_cats}We will be interested in metric spaces which are also commutative monoids for which the metric and monoid structures are in some sense compatible. This idea is formalized categorically by the notion of a commutative monoid internal to a category. In order to make this notion precise, the category in question must have additional structure - that of a symmetric monoidal category.
  For example, a commutative topological monoid, analogous to a topological group, is a 
  commutative monoid
internal to $\cat{Top}$, the (symmetric monoidal) category of topological spaces and continuous maps. 
In this section we sketch the formal definition of symmetric monoidal categories and commutative monoid objects internal to them. For a complete treatment, see \cite[VII.1]{mac1998categories}.

A \emph{symmetric monoidal category} is a category $\cat{C}$ equipped with a functor $\otimes:\cat{C}\times \cat{C}\to \cat{C}$ called the \emph{tensor product}, an object $1\in \cat{C}$ called the \emph{unit object}, a natural isomorphism $\alpha_{X,Y,Z}:(X\otimes Y)\otimes Z\to X\otimes (Y\otimes Z)$ called the \emph{associator}, a natural isomorphism $\lambda_X:1\otimes X\to X$ called the \emph{left unitor}, a natural isomorphism $\rho_X:X\otimes 1\to X$ called the \emph{right unitor}, and a natural isomorphism $B_{X,Y}:X\otimes Y\to Y\otimes X$ called the \emph{braiding}. These natural isomorphisms must satisfy certain \emph{coherence conditions} expressed by commutative diagrams. For example, the unitors and the associator must obey the \emph{triangle equality}, which specifies that the diagram 
\[
	\begin{tikzcd}
		(X\otimes 1)\otimes Y \ar[rr,"\alpha_{X,1,Y}"] \ar[dr,swap,"\rho_X\otimes \textup{id}_Y"]& &X\otimes (1\otimes Y)\ar[dl,"\textup{id}_X \otimes \lambda_X"] \\
		& X\times Y
	\end{tikzcd}
\]
commutes. For the other coherence conditions and their diagrams, see \cite[VII.1]{mac1998categories}.

\begin{example}\label{ex:set_is_sym_mon_cat}
	$\cat{Set}$ is a symmetric monoidal category with tensor product given by the cartesian product and with unit object being the one-point set $*$. The associator, left unitor, right unitor, and braiding are defined by the obvious bijections $(X\times Y)\times Z \cong X\times (Y\times Z)$, $*\times X\cong X$, $X\times *\cong X$, and $X\times Y\cong Y\times X$, respectively.
\end{example}

\begin{example}\label{ex:met_is_sym_mon_cat}
	$\cat{Met}$ is a symmetric monoidal category with tensor product $(X,d_X)\otimes (Y,d_Y)\allowbreak := (X\times Y, D_\infty)$, where $D_\infty((x,y),(x',y')):= \max(d_X(x,x'),d_Y(y,y'))$, and with unit object the one-point metric space $*$. The associator, unitors, and braiding are the same as those in $\cat{Set}$ (it only needs to be checked that these maps are isometries, i.e., isomorphisms in $\met$).
\end{example}

For any symmetric monoidal category $(\cat{C},\otimes,1)$, a notion of a commutative monoid defined within $\cat{C}$ can be made precise. A \emph{commutative monoid object} in $(\cat{C},\otimes, 1)$ (or a \emph{commutative monoid internal to $(\cat{C},\otimes,1))$} is a tuple $(M,\mu,e)$ consisting of an object $M\in \cat{C}$, a morphism $e:1\to M$ called the \emph{unit}, and a morphism $\mu:M\otimes M\to M$ called the \emph{product}. The unit and product morphisms are required to satisfy certain coherence conditions with the associator, unitors, and braiding of $\cat{C}$, expressing associativity and commutativity of the product and the fact that the unit morphism serves as an identity for the product. For example, the unit $e:1\to M$ is required to make the diagram 
\[
\begin{tikzcd}
	1\otimes M\ar[r,"e\otimes \textup{id}_M"] \ar[dr,swap,"\lambda_M "]& M\otimes M \ar[d,"\mu"] \ar[r,"\textup{id}_M\otimes e"]& M\otimes 1 \ar[dl,"\rho_M"] \\
	& M
\end{tikzcd}
\]
commute. See again \cite[VII.1]{mac1998categories} for the other coherence conditions.

A morphism between commutative monoid objects $(M,\mu,e)$, $(M',\mu',e')$ in a symmetric monoidal category $(\cat{C},\otimes,1)$ is a morphism $f:M\to M'$ in $\cat{C}$ such that $f\mu = \mu'(f\otimes f)$ and $fe = e'$. Commutative monoid objects in $(\cat{C}, \otimes, 1)$, together with these morphisms, form a category which is denoted by $\cat{CMon}(\cat{C},\otimes,1)$ or more simply by $\cat{CMon(C)}$ when the tensor product and unit are fixed.

\begin{example}
	Recall that $\cat{Set}$ is a symmetric monoidal category with tensor product the cartesian product and with unit object the one-point set (see Example \ref{ex:set_is_sym_mon_cat}). The commutative monoid objects internal to $\cat{Set}$ are just commutative monoids, and the corresponding morphisms are monoid homomorphisms. That is, we have an isomorphism of categories $\cat{CMon(Set)}\cong \cat{CMon}$.
\end{example}

\begin{example}
	Let $\met$ be given the symmetric monoidal structure of Example \ref{ex:met_is_sym_mon_cat}. A commutative monoid internal to $\met$ is a metric space $(M,d)$ together with a unit $e:*\to \met$ and a multiplication operation $+:M\times M\to M$. Note that $e$ simply picks out an element of $M$, which we will denote by $0$. The fact that $+$ is metric map means that 
	\begin{equation}\label{eq:met_mon}d(a+b,a'+b')\leq \max(d(a,a'),d(b,b')) \quad \textup{ for all $a,a',b,b'\in M$}.
	\end{equation}Thus a commutative monoid internal to $\met$ is a tuple $(M,d,+,0)$, where $(M,d)$ is a metric space, $(M,+,0)$ is a commutative monoid, and such that \eqref{eq:met_mon} holds.
\end{example}

We will consider variations of the preceding example for pointed metric spaces and for different choices of tensor product.

\section{Metric pairs and pointed metric spaces}
\label{sec:metric-spaces}
In this section , we introduce the categories of metric spaces of interest to us. In Section \ref{sec:metric_pairs}, we introduce the categories of metric pairs and pointed metric spaces. We show that every metric pair gives rise to a pointed metric space by taking a quotient and that this construction is functorial. Moreover, we show that the category of pointed metric spaces is a reflective subcategory of the category of metric pairs, with reflector being the quotient functor. In Section \ref{sec:p-strengthening}, we introduce the $p$-strengthened triangle inequality ($p\in [1,\infty])$ and the corresponding subcategories of pointed metric spaces which satisfy it. These subcategories are required for our statement of universality.

\subsection{Metric pairs, pointed metric spaces, and quotients}\label{sec:metric_pairs}
We now introduce the main categories of interest to us.
\begin{definition} \label{def:metric-pair}
  Let $\cat{Met_{pairs}}$ denote the category whose objects are of the form $(X,d,A)$, where $(X,d)$ is a metric space with $A$ a nonempty subset of $X$ and whose morphisms $f:(X,d,A) \to (Y,d',B)$ are metric maps $f:(X,d) \to (Y,d')$ such that $f(A) \subset B$. $(X,d,A)$ is called a \emph{metric pair}.
\end{definition}

\begin{example} \label{ex:half-plane}
  Consider the metric space $(\R^2,d)$ where $d$ is the metric induced by the $q$-norm, where $1 \leq q \leq \infty$. Let $\R^2_{\leq} = \{(x,y) \in \R^2 \st x \leq y\}$ and similarly define subsets $\R^2_{\geq}$ and $\R^2_{=}$ where the latter is also denoted $\Delta$. Then we have metric pairs
  $(\R^2,d,\R^2_{\geq})$ and $(\R^2_{\leq},d,\Delta)$.
  Similarly, we have the metric pair
  $(\eR^2_{\leq},d,\Delta)$.
\end{example}

\begin{definition} \label{def:pointed-metric-space}
  Let $\cat{Met_*}$ denote the full subcategory of $\cat{Met_{pairs}}$ whose objects are of the form $(X,d,\{x_0\})$, which we denote $(X,d,x_0)$. We call $x_0$ the \emph{basepoint} and call $(X,d,x_0)$ a \emph{pointed metric space}.
  A morphism $f:(X,d,x_0) \to (Y,d',y_0)$ is called a \emph{pointed metric map}.
\end{definition}

\begin{example} \label{ex:intervals}
  Let $\Int(\R)$ denote the set of intervals in $\R$ with $d(I,J)$ equal to the length (i.e. Lebesgue measure) of the symmetric difference
  $(I \cup J) \setminus (I \cap J)$. Then $(\Int(\R),d,\emptyset)$ is a pointed metric space.
  We may also equip $\Int(\R)$ with the {Hausdorff distance} $d_H$
  to obtain the pointed metric space $(\Int(\R),d_H,\emptyset)$.
\end{example}

We now show how to obtain a pointed metric space from a metric pair.

\begin{definition} \label{def:quotient-pointed-metric-space}
  Given a metric pair $(X,d,A)$ consider the quotient set $X/A = (X \setminus A) \amalg \{A\}$.
  Let $\overline{d}:X/A \times X/A \to [0,\infty]$ be the induced metric. That is, $\overline{d}(A,A) = 0$,
  for $x \in X \setminus A$, $\overline{d}(x,A) = \overline{d}(A,x) = d(x,A)$, where $d(x,A) = \inf_{y \in A}d(x,y)$,
  and for $x,y \in X \setminus A$, $\overline{d}(x,y) = \min\left(d(x,y),d(x,A) + d(y,A)\right)$.
  
  Given a morphism $f:(X,d_X,A) \to (Y,d_Y,B)$, let $\overline{f}:X/A \to Y/B$ be the induced map. That is, $\overline{f}(A) = B$, and for $x \in X \setminus A$, $\overline{f}(x) = B$ if $f(x) \in B$ and otherwise $\overline{f}(x) = f(x)$. 
  
  There is a natural quotient map of pairs $q:(X,A)\to (X/A,\{A\})$ given by $q(x) = x$ if $x\in X\setminus A$ and $q(x) = \{A\}$ if $x\in A$.  We will sometimes denote the image of $x$ by $[x]$, but we will drop the brackets when there is no ambiguity.
\end{definition}

We will show that this quotient map may be used to define a functor from $\cat{Met_{pairs}}$ to $\cat{Met_*}$. First, we show that the quotient of a metric pair is indeed a pointed metric space.

\begin{lemma} 
  If $(X,d,A)$ is a metric pair then $(X/A,\overline{d},A)$ is a pointed metric space. Moreover, the quotient map $q:(X,d,A)\to (X/A,\overline{d},A)$ is a metric map.
\end{lemma}

\begin{proof}
  We need to show that $\overline{d}$ is a metric. Point triviality and symmetry follow from the definition. It remains to prove the triangle inequality.
  There are three nontrivial cases.
  
  For the first case, let $x,z \in X \setminus A$. We want to show that $\overline{d}(x,A) \leq \overline{d}(x,z) + \overline{d}(z,A)$. Since $\overline{d}(x,z) + \overline{d}(z,A) = \min(d(x,z), d(x,A) + d(z,A)) + d(z,A)$, it suffices to show that $\overline{d}(x,A)\leq d(x,z) + d(z,A)$ and $\overline{d}(x,A) \leq d(x,A) + 2d(z,A)$.  
  The first inequality holds since $\overline{d}(x,A) = d(x,A) = \inf_{y\in A} d(x,y) \leq \inf_{y\in A}(d(x,z) + d(z,y)) = d(x,z) + d(z,A)$, and the second inequality holds trivially since
  $\overline{d}(x,A) = d(x,A) \leq d(x,A) + 2d(z,A)$.
  Thus $\overline{d}(x,A) \leq \overline{d}(x,z) + \overline{d}(z,A)$.

  For the second case, let $x,y \in X \setminus A$.
  Then $\overline{d}(x,y) = \min( d(x,y), d(x,A)+d(A,y) )
  \leq d(x,A) + d(A,y) = \overline{d}(x,A) + \overline{d}(A,y)$.

  For the third case, let $x,y,z \in X \setminus A$.
  We want to show that
  $\overline{d}(x,y) \leq \overline{d}(x,z) + \overline{d}(z,y)
  = \min( d(x,z), d(x,A) + d(A,z) ) + \min( d(z,y), d(z,A) + d(A,y) )$.
  The right hand side has four possible values. First,
  $\overline{d}(x,y) \leq d(x,y) \leq d(x,z) + d(z,y)$. Second,
  $\overline{d}(x,y) \leq d(x,A) + d(A,y)
  \leq d(x,z) + d(z,A) + d(A,y)$. Third,
  $\overline{d}(x,y) \leq d(x,A) + d(A,y)
  \leq d(x,A) + d(A,z) + d(z,y)$. Fourth,
  $\overline{d}(x,y) \leq d(x,A) + d(A,y)
  \leq d(x,A) + d(A,z) + d(z,A) + d(A,y)$.
  Therefore
    $\overline{d}(x,y) \leq \overline{d}(x,z) + \overline{d}(z,y)$.
    
   To prove the second statement, let $x,y\in X$. There are three cases. First, if $x,y\in X\setminus A$ then $\overline{d}(q(x),q(y)) = \min(d(x,y),d(x,A)+d(y,A)) \leq d(x,y)$. Second, if $x\in X\setminus A$ and $y\in A$, then $\overline{d}(q(x),q(y)) = d(x,A)\leq d(x,y)$. Third, if $x,y\in A$ then $\overline{d}(q(x),q(y)) = \overline{d}(A,A) = 0 \leq d(x,y)$. This completes the proof.
   \end{proof}

Next, we show that this map sends morphisms in $\cat{Met_{pairs}}$ to morphisms in $\cat{Met_*}$.
\begin{lemma}
  Given a morphism $f:(X,d_X,A) \to (Y,d_Y,B)$ of metric pairs, the induced map $\overline{f}:(X/A,\overline{d}_X,A) \to (Y/B,\overline{d}_Y,B)$ is a pointed metric map.
\end{lemma}

\begin{proof}
  We will prove that $\overline{f}$ is a metric map.
  Let $x \in X \setminus A$.
  First we show that $d_Y(f(x),B) \leq d_X(x,A)$.
  Indeed, $d_Y(f(x),B) = \inf_{y \in B} d_Y(f(x),y) \leq \inf_{x' \in A} d_Y(f(x),f(x')) \leq \inf_{x' \in A} d_X(x,x') = d(x,A)$.
  Then $\overline{d}_Y(\overline{f}(x), \overline{f}(A)) = \overline{d}_Y(f(x),B) = d_Y(f(x),B) \leq d_X(x,A) = \overline{d}_X(x,A)$.
  Next let $x,x' \in X \setminus A$.
  Then $\overline{d}_Y(\overline{f}(x),\overline{f}(x'))
  = \overline{d}_Y(f(x),f(x'))
  = \min( d_Y(f(x),f(x')), d_Y(f(x),B) + d_Y(B,f(x')) )
  \leq \min( d_X(x,x'), d_X(x,A) + d_X(A,x') ) = \overline{d}_X(x,x')$.
\end{proof}

With the above lemmas in hand, it is now easy to check that we have a functor.

\begin{definition} \label{def:Q}
  Let $Q: \cat{Met_{pairs}} \to \cat{Met_*}$ be the functor that sends a metric pair $(X,d,A)$ to the pointed metric space $(X/A,\overline{d},A)$ and that sends
   $f:(X,d_X,A) \to (Y,d_Y,B)$ to $\overline{f}:(X/A,\overline{d}_X,A) \to (Y/B,\overline{d}_Y,B)$.
\end{definition}

\begin{theorem} \label{thm:Q}
  $\cat{Met_*}$ is a reflective subcategory of $\cat{Met_{pairs}}$ with left adjoint $Q$.
\end{theorem}

\begin{proof}

  By Lemma \ref{lem:adjoint}, it suffices to show that, for each $X = (X,d_X,A)\in \cat{Met_{pairs}}$, there is a morphism $r:(X,d_X,A)\to (X/A,\overline{d},A) = IQ(X,d_X,A)$ in $\cat{Met_{pairs}}$ for which $((X/A,\overline{d},A),r)$ is an initial object in the comma category $X\downarrow I$. Here, $I:\cat{Met_*}\to\cat{Met_{pairs}}$ denotes the inclusion functor.
  \begin{equation} \label{cd:1}
  	\begin{tikzcd}
  		(X,d_X,A) \ar[r,"r"] \ar[dr,swap,"f"]& (X/A,\overline{d}_X,A) \ar[d,dashed,"\exists ! \overline{f}"]\\
  		& (Y,d_Y,y_0)
  	\end{tikzcd}
  \end{equation}
   To this end, define $r:(X,d_X,A) \to (X/A,\overline{d}_X,A)$ by $r(x) = A$ if $x \in A$ and $r(x) = x$ otherwise.
  Let us show that $r$ is a metric map.
  Let $x,x' \in X$.
  If $x,x' \not\in A$ then $\overline{d}_X(rx,rx') = \overline{d}_X(x,x') = \min(d_X(x,x'),d_X(x,A)+d_X(A,x)) \leq d_X(x,x')$.
  If $x \not\in A$ and $x' \in A$ then $\overline{d}_X(rx,rx') = \overline{d}_X(x,A) = d_X(x,A) \leq d_X(x,x')$.

  To see that $r$ is a universal element, let $f:(X,d_X,A) \to (Y,d_Y,y_0) \in \cat{Met_{pairs}}$ be given. We want to show that there is a unique metric map $\overline{f}:(X/A,\overline{d}_X,A)\to (Y,d_y,y_0)$ such that $\overline{f}\circ r = f$. By the commutativity of \eqref{cd:1}, we are forced to define $\overline{f}(A) = y_0$ and $\overline{f}(x) = f(x)$ if $x \not\in A$. This demonstrates uniqueness.
  To establish existence, it remains to show that $\overline{f}$ is a metric map.
  If $x \in X \setminus A$ then $d_Y(\overline{f}(x),\overline{f}(A)) = d_Y(f(x),y_0) = d_Y(f(x),f(x')) \leq d_X(x,x')$ for all $x' \in A$.
  Thus $d_Y(\overline{f}(x),\overline{f}(A)) \leq d_X(x,A) = \overline{d}_X(x,A)$.
  If $x,x' \in X \setminus A$ then
  $d_Y(\overline{f}(x),\overline{f}(x')) = d_Y(f(x),f(x')) \leq d_X(x,x')$.
  Furthermore $d_Y(\overline{f}(x),\overline{f}(x')) = d_Y(f(x),f(x')) \leq d_Y(f(x),y_0)+d_Y(y_0,f(x')) \leq d_X(x,A) + d_X(A,x')$.
  Therefore $d_Y(\overline{f}(x),\overline{f}(x')) \leq \overline{d}_X(x,x')$.
\end{proof}

\subsection{The $p$-strengthened triangle inequality}\label{sec:p-strengthening}
In this section, we introduce a convenient class of pointed metric spaces for each $p\in [1,\infty]$. These are pointed metric spaces which satisfy a slightly stronger version of the triangle inequality with respect the basepoint.

\begin{definition}\label{def:p-strengthened}
  Let $(X,d)$ be a metric space, $x_0\in X$, and $p\in [1,\infty]$. We say that the metric $d$ satisfies the \emph{$p$-strengthened triangle inequality with respect to $x_0$} if $d(x,y)\leq \elp{(d(x,x_0),d(x_0,y))}{p}$ for all $x,y\in X$.
  Let $\cat{Met_*^p}$ denote the full subcategory of $\cat{Met_*}$ consisting of those objects $(X,d,x_0)$ for which $(X,d)$ satisfies the $p$-strengthened triangle inequality with respect to $x_0$.
\end{definition}

Note that the $1$-strengthened triangle inequality is just the triangle inequality. 
So $\cat{Met_*^1} = \cat{Met_*}$.
Also, for $1 \leq p \leq q \leq \infty$, $\cat{Met_*^q}$ is a full subcategory of $\cat{Met_*^p}$.

\begin{example}
  Let $*$ denote the singleton set. Also let $*$ denote the pointed metric space $(*,0,*)$. Then for all $p \in [1,\infty]$, $* \in \cat{Met_*^p}$. In fact, it is the initial and terminal object in $\cat{Met_*^p}$.
\end{example}

From any metric space $(X,d)$ and given basepoint $x_0\in X$, we can obtain a metric that satisfies the $p$-strengthened inequality. The following definition is a more general construction.

\begin{definition} \label{def:dp}
  Let $(X,d,A)$ be a metric pair and let $p \in [1,\infty]$.
  Define $d_p:X \times X \to [0,\infty]$ by
  \[d_p(x,y) := \min(d(x,y),\elp{(d(x,A),d(A,y))}{p}).\]
\end{definition}

In the special case that $A = \{x_0\}$ is a singleton, it is clear from the definition that $d_p$ satisfies the $p$-strengthened triangle inequality with respect to $x_0$. We still need to verify that $d_p$ is actually a metric.

\begin{lemma}
  Let $p \in [1,\infty]$.
  If $(X,d,A)$ is a metric pair then so is $(X,d_p,A)$.
\end{lemma}

\begin{proof}
Point triviality and symmetry follow from the definition.
  To show the triangle inequality, let $x,y,z \in X$.
  We want to show that
  $\min( d(x,y), \norm{ ( d(x,A), d(A,y) ) }_p )
  \leq \min( d(x,z), \norm{ ( d(x,A), d(A,z) ) }_p )
  + \min( d(z,y), \norm{ ( d(z,A), d(A,y) ) }_p )$.
  The right hand side has four possible values. 
  First, $d(x,y) \leq d(x,z) + d(z,y)$.
  Second,
\begin{multline*}\norm{ (d(x,A),d(A,y)) }_p
  \leq \norm{ (d(x,z)+d(z,A),d(A,y)) }_p \\
  = \norm{ (d(x,z),0) + (d(z,A),d(A,y)) }_p 
  \leq d(x,z) + \norm{ (d(z,A),d(A,y)) }_p.
  \end{multline*}
  Third, \begin{multline*}\norm{ (d(x,A),d(A,y)) }_p
  \leq \norm{ (d(x,A),d(A,z)+d(z,y)) }_p\\
  = \norm{ (d(x,A),d(A,z)) + (0,d(z,y)) }_p 
  \leq \norm{ (d(x,A),d(A,z)) }_p + d(z,y).\end{multline*}
  Fourth, \begin{multline*}\norm{ (d(x,A),d(A,y)) }_p
  \leq \norm{ (d(x,A),d(A,z),d(z,A),d(A,y)) }_p\\
  = \norm{ ( \norm{(d(x,A),d(A,z))}_p,\norm{(d(z,A),d(A,y))}_p) }_p\\
  \leq \norm{(d(x,A),d(A,z))}_p +\norm{(d(z,A),d(A,y))}_p.
\end{multline*}
The result now follows from these inequalities.
\end{proof}

\begin{lemma}
  Let $p \in [1,\infty]$.
  If $f:(X,d,x_0) \to (Y,d',y_0)$ is a pointed metric map, then so is $f:(X,d_p,x_0) \to (Y,d'_p,y_0)$.
\end{lemma}

\begin{proof}
  Let $x,x' \in X$.
  Then \begin{multline*}
  	d'_p(f(x),f(x'))
  = \min( d'(f(x),f(x')), \norm{ ( d'(f(x),y_0), d'(y_0,f(x')) ) }_p )\\
  = \min( d'(f(x),f(x')), \norm{ ( d'(f(x),f(x_0)), d'(f(x_0),f(x')) ) }_p )\\
  \leq \min( d(x,x'), \norm{ ( d(x,x_0), d(x_0,x') ) }_p )
  = d_p(x,x').
    \end{multline*}
\end{proof}

The operation which sends a pointed metric space $(X,d,x_0)$ to the $(X,d_p,x_0)$ is easily seen to be functorial.
\begin{definition} \label{def:Spq}
  Let $1 \leq p \leq q \leq \infty$.
  Let $S_{p,q}: \cat{Met_*^p} \to \cat{Met_*^q}$ be the functor that sends $(X,d,x_0)$ to $(X,d_q,x_0)$ and $f:(X,d,x_0) \to (Y,d',y_0)$ to $f:(X,d_q,x_0) \to (Y,d'_q,y_0)$.  We will also denote $S_{1,p}$ by $S_p$.
\end{definition}

\begin{theorem} \label{thm:Spq}
  Let $1 \leq p \leq q \leq \infty$.
  $\cat{Met_*^q}$ is a reflective subcategory of $\cat{Met_*^p}$ with left adjoint $S_{p,q}$.
  As a special case, 
  $\cat{Met_*^p}$ is a reflective subcategory of $\cat{Met_*}$ with left adjoint $S_{p}$.
\end{theorem}

\begin{proof}
  Let $(X,d,x_0) \in \cat{Met_*^p}$, $(Y,d',y_0) \in \cat{Met_*^q}$ and $f:(X,d,x_0) \to (Y,d',y_0) \in \cat{Met_*^p}$.
  \begin{equation} \label{cd:2}
    \begin{tikzcd}
      (X,d,x_0) \ar[r,"r"] \ar[dr,swap,"f"]& (X,d_q,x_0) \ar[d,dashed,"\exists ! \overline{f}"]\\
      & (Y,d',y_0)
    \end{tikzcd}
  \end{equation}

    Let $r:(X,d,x_0) \to (X,d_q,x_0) \in \cat{Met_*^p}$ be the identity function on $X$. We will show that $((X,d_q,x_0),r)$ is a universal element in the comma category $X\downarrow I$, where $I$ denotes the inclusion $\cat{Met_*^q}\hookrightarrow \cat{Met_*^p}$. Note that $r$ is a metric map since for all $x,x' \in X$,
  $d_q(x,x') \leq d(x,x')$.
  Let $f:(X,d,A)\to (Y,d',y_0)$ be a pointed metric map. By the commutativity of \eqref{cd:2}, we are forced to define $\overline{f} = f$. This establishes uniqueness. To establish existence, it remains to show that $\overline{f}$ is a metric map.
  Let $x,x' \in X$.
  We have $d'(f(x),f(x')) \leq d(x,x')$ and, since $(Y,d',y_0)\in \cat{Met_*^q}$,

\begin{multline*}d'(f(x),f(x'))
  \leq \norm{ ( d'(f(x),y_0), d'(y_0,f(x')) ) }_q\\
  = \norm{ ( d'(f(x),f(x_0)), d'(f(x_0),f(x')) ) }_q
  \leq \norm{ ( d(x,x_0), d(x_0,x') ) }_q.
  \end{multline*}
  Therefore $d'(\overline{f}(x),\overline{f}(x')) \leq d_q(x,x')$. It follows by Lemma \ref{lem:adjoint} that $S_{p,q}$ is left adjoint to the inclusion.
\end{proof}

Combining Definitions \ref{def:quotient-pointed-metric-space} and \ref{def:dp}, it is easy to verify the following.

\begin{lemma} \label{lem:overlineDp}
  Given a metric $(X,d,A)$ and $p \in [1,\infty]$, we have the pointed metric space $(X/A,\overline{d}_p,A)$ where $\overline{d}_p$ is given by $\overline{d}_p(A,A) = 0$, for $x \in X \setminus A$, $\overline{d}_p(x,A) = d(x,A) = \overline{d}_p(A,x)$, and for $x,y \in X \setminus A$, $\overline{d}_p(x,y) = \min( d(x,y), \norm{d(x,A),d(A,y)}_p)$.
\end{lemma}

\section{Commutative metric monoids and Wasserstein distance}
\label{sec:metric-monoids-wasserstein}

In this section, we introduce $p$-subadditive commutative metric monoids. These are metric spaces which are also monoids and for which the monoid operation is, in a precise sense, compatible with the metric. The category of $p$-subadditive commutative metric monoids is the setting in which we state our universality results for the Wasserstein distances.
	

\subsection{Commutative metric monoids}
\label{sec:metric-monoids} In this section, we introduce $p$-subadditive commutative metric monoids which are metric spaces that are simultaneously commutative monoids and which satisfy a certain compatibility condition between the metric and monoid operation.

\begin{definition} \label{def:metric-monoid}

    Let $p \in [1,\infty]$.
  A \emph{$p$-subadditive commutative metric monoid} is a tuple $(M,d,+,0)$ where $(M,d,0)$ is a pointed metric space and $(M,+,0)$ is a commutative monoid such that for all $a,b,a',b' \in M$,

\begin{equation*}
    d(a+b,a'+b') \leq \norm{(d(a,a'),d(b,b'))}_p.
  \end{equation*}
  In this case, we say that the metric $d$ is \emph{$p$-subadditive}. A \emph{morphism of $p$-subadditive commutative metric monoids} $f:(M,d,+,0) \to (N,\rho,+,0)$ is a pointed metric map $f:(M,d,0) \to (N,\rho,0)$ such that $f:(M,+,0) \to (N,+,0)$ is a monoid homomorphism.
  Call such a map a \emph{metric monoid homomorphism}.
  Let
  $\cat{CMetMon_p}$
  denote the category of $p$-subadditive commutative metric monoids and metric monoid homomorphisms.
\end{definition}


\begin{lemma}
  Let $p \in [1,\infty]$.
  Let $(M,d,+,0)$ be a $p$-subadditive commutative metric monoid.
  Then for $n \geq 0$ and $a_1,\ldots,a_n,b_1,\ldots,b_n \in M$,
  \begin{equation*}
    d(a_1+\cdots+a_n,b_1+\cdots+b_n) \leq \norm{(d(a_i,b_i))_{i=1}^n}_p.
  \end{equation*}
\end{lemma}

\begin{proof}
  The proof is by induction on $n$.
  \begin{multline*}
    d(a_1 + \cdots + a_{n+1},b_1+\cdots+b_{n+1}) 
  \leq \norm{ (d(a_1 + \cdots + a_n,b_1+\cdots+b_n), d(a_{n+1},b_{n+1}) )}_p\\
  \leq \norm{ (\norm{(d(a_i,b_i))_{i=1}^n}_p, d(a_{n+1},b_{n+1}) )}_p
  = \norm{(d(a_i,b_i))_{i=1}^{n+1}}_p.\qedhere
  \end{multline*}
\end{proof}

\begin{corollary} \label{cor:p-subadditive}
  Let $p \in [1,\infty]$.
  Let $(M,d,+,0)$ be a $p$-subadditive commutative metric monoid.
  Then for $n \geq 0$ and $a_1,\ldots,a_n,b_1,\ldots,b_n \in M$,
  \begin{equation*}
    d(a_1+\cdots+a_n,b_1+\cdots+b_n) \leq \min_{\sigma \in \Sigma_n} \norm{(d(a_i,b_{\sigma(i)}))_{i=1}^n}_p,
  \end{equation*}
	where $\Sigma_n$ denotes the symmetric group on $n$ symbols.
\end{corollary}

The following lemma shows that there is a forgetful functor 
$U_p:\cat{CMetMon_p} \to \cat{Met_*^p}$.

\begin{lemma}\label{lem:forgetful_functor_U}
  Let $p \in [1,\infty]$.
  Let $(M,d,+,0)$ be a $p$-subadditive commutative metric monoid. Then $(M,d,0) \in \cat{Met_*^p}$.
  Furthermore if $f:(M,d,+,0) \to (N,\rho,+,0)$ is a morphism of commutative metric monoids then $f:(M,d,0) \to (N,\rho,0) \in \cat{Met_*^p}$.
\end{lemma}

\begin{proof}
  For the first statement, we show that $(M,d,0)$ satisfies the $p$-strengthened triangle inequality at $0$. Let $a,b \in M$. Then
  $d(a,b) = d(a+0,0+b) \leq \norm{(d(a,0),d(0,b))}_p$.
  The second statement follows directly from the definitions.
\end{proof}

\subsection{Monoid objects in $\cat{Met_*^p}$}
\label{sec:monoidal-cat}
In this section, we show that $\cat{Met_*^p}$ can be equipped with a tensor product making it into a symmetric monoidal category. We then show that $p$-subadditive commutative metric monoids are precisely the commutative monoids internal to this symmetric monoidal category.
\begin{definition} \label{def:xp}
  Let $p \in [1,\infty]$.
  Given pointed metric spaces $(X,d_X,x_0)$ and $(Y,d_Y,y_0)$, define $d_X \times_p d_Y: (X \times Y) \times (X \times Y) \to [0,\infty]$ by
  \begin{equation*}
    (d_X \times_p d_Y) \left( (x,y), (x',y') \right) = \norm{(d_X(x,x'),d_Y(y,y'))}_p.
  \end{equation*}
  We call $d_X\times_p d_Y$ the \emph{$p$-product metric}. Let $X \times_p Y$ denote the tuple $(X\times Y, d_X \times_p d_Y, (x_0,y_0))$.
\end{definition}
The following lemma shows that $d_X\times_p d_Y$ defines a metric on the product $X\times Y$.
\begin{lemma}
  Let $p \in [1,\infty]$.
  If $X = (X,d_X,x_0)$ and $Y = (Y,d_Y,y_0)$ are pointed metric spaces then so is $X \times_p Y$.
\end{lemma}

\begin{proof}

    We show that $d_X \times_p d_Y$ is a metric for $X\times Y$.
  Point triviality and symmetry follow from the corresponding properties for $d_X$ and $d_Y$.
  It remains to prove the triangle inequality.
  For all $(x,y), (x',y'), (x'',y'')\in X \times Y$, 

\begin{align*}
    \norm{(d_X(x,x'),d_Y(y,y'))}_p
    &\leq \norm{(d_X(x,x'')+d_X(x'',x'),d_Y(y,y'')+d_Y(y'',y'))}_p\\
    &= \norm{(d_X(x,x''),d_Y(y,y'')) + (d_X(x'',x'),d_Y(y'',y'))}_p\\
    &\leq \norm{(d_X(x,x''),d_Y(y,y''))}_p + \norm{(d_X(x'',x'),d_Y(y'',y'))}_p. \qedhere
  \end{align*}
\end{proof}

The product metric $d_X\times_p d_Y$ can be used to give a succinct and convenient description of $p$-subadditive commutative metric monoids.
\begin{lemma}
  A \emph{$p$-subadditive commutative metric monoid} is a tuple $(M,d,+,0)$ where $(M,d,0)$ is a pointed metric space and $(M,+,0)$ is a commutative monoid such that $+:M \times_p M \to M$ is a metric map.
\end{lemma}

\begin{proof}
  Consider $(M,d,+,0)$ where $(M,d,0)$ is a pointed metric space and $(M,+,0)$ is a commutative monoid.
  $+:M \times_p M \to M$ is a metric map 
  if and only if
  for all $(a,b),(a',b') \in M \times M$,
  $d(a+b,a'+b') \leq \norm{d(a,a'),d(b,b')}_p$.
\end{proof}

The following lemma shows that $\cat{Met_*^p}$ is closed with respect to forming $p$-product metrics.
\begin{lemma}\label{lem:met_p_closed}
  Let $p \in [1,\infty]$.
  If $X = (X,d_X,x_0), Y = (Y,d_Y,y_0) \in \cat{Met_*^p}$ then $X \times_p Y \in \cat{Met_*^p}$.
\end{lemma}

\begin{proof}
  We show that $d_X \times_p d_Y$ satisfies the $p$-strengthened triangle inequality with respect to $(x_0,y_0)$.
  For all $(x,y), (x',y') \in X \times Y$, 
  \begin{align*}
    \norm{(d_X(x,x'),d_Y(y,y'))}_p
    &\leq \norm{ ( \norm{(d_X(x,x_0),d_X(x_0,x'))}_p, \norm{(d_Y(y,y_0),d_Y(y_0,y'))}_p ) }_p\\
    &= \norm{ (d_X(x,x_0),d_X(x_0,x'),d_Y(y,y_0),d_Y(y_0,y')) }_p\\
    &= \norm{ ( \norm{(d_X(x,x_0),d_Y(y,y_0))}_p, \norm{(d_X(x_0,x'),d_Y(y_0,y'))}_p ) }_p.\qedhere
  \end{align*}
\end{proof}

We want to show that $(\cat{Met_*^p}, \times_p, *)$ is a symmetric monoidal category. We will first show that $(\cat{Met_*}, \times_p, *)$
  is a symmetric monoidal category.

\begin{proposition} \label{prop:monoidal-category}
		For each $p\in [1,\infty]$, $(\cat{Met_*},\times_p,*)$ is a symmetric monoidal category.
\end{proposition}
	
\begin{proof}
  For the associator, consider $X,Y,Z \in \cat{Met_*}$ and
  $x,x' \in X$, $y,y' \in Y$ and $z,z' \in Z$. Then
  \begin{multline*}
    ((d_X \times_p d_Y) \times_p d_Z)(((x,y),z),((x',y'),z))
    = \norm{((d_X \times_p d_Y)((x,y),(x',y')),d_Z(z,z'))}_p\\
    = \norm{( \norm{(d_X(x,x'),d_Y(y,y'))}_p, d_Z(z,z'))}_p = \norm{(
      d_X(x,x'), d_Y(y,y'), d_Z(z,z') )}_p.
  \end{multline*}
  The left unitor is an isometry since for $X \in \cat{Met_*}$ and $x,x' \in X$, $\norm{(d_X(x,x'),0)}_p = d_X(x,x')$.
  Similarly, the right unitor is an isometry. The braiding is given by
    the obvious isometry $X\times_p Y \cong Y\times_p X$.
  With these computations in
  hand, the rest of the axioms are easy to check.
\end{proof}
	
\begin{corollary} \label{cor:monoidal-category} For each
  $p\in [1,\infty]$, $(\cat{Met_*^p},\times_p,*)$ is a symmetric
  monoidal category, which we denote by $\cat{Met_*^p}$.
		
  \begin{proof} This follows immediately since $\cat{Met_*^p}$ is
      a subcategory of $\cat{Met_*}$ which, by Lemma
      \ref{lem:met_p_closed}, is closed under the tensor product
      $\times_p$.
    \end{proof}
\end{corollary}

Let $\cat{CMon(Met_*,\times_p,*)}$
and $\cat{CMon(Met_*^p)}$ denote the categories of
commutative monoids internal to the symmetric monoidal categories
$(\cat{Met_*},\times_p,*)$
and $\cat{Met_*^p}$, respectively.
Recall that $\cat{CMetMon_p}$ denotes the category of
  $p$-subadditive commutative metric monoids. The following
  proposition shows that $\cat{CMetMon_p}$ is precisely the
  category of commutative monoids internal to the $\cat{Met_*^p}$.
  Moreover, we show that
  $\cat{CMon(Met_*,\times_p,*)}$ and $\cat{CMon(Met_*^p)}$ are in fact
  the same.

\begin{theorem} \label{thm:CMon}
  A commutative monoid in the symmetric monoidal category
  $\cat{Met_*^p}$ is a $p$-subadditive commutative metric monoid and a
  morphism of commutative monoids in $\cat{Met_*^p}$ is morphism of
  $p$-subadditive commutative metric monoids.
  That is, $\cat{CMetMon_p} = \cat{CMon(Met_*^p)}$.
  Moreover,
  $\cat{CMon(Met_*,\times_p,*)} = \cat{CMon(Met_*^p)}$.
\end{theorem}
              
\begin{proof}
  A commutative monoid in 
  $\cat{Met_*^p}$ is a pointed metric space
  $(M,d,m_0) \in \cat{Met_*^p}$ together with a binary operation
  $+: M \times M \to M$ that is associative, commutative, and for
  which $m_0$ is a unit, such that
  $+: M \times_p M \to M \in \cat{Met_*^p}$.  That is, $+$ is
  $p$-subadditive. Thus a commutative monoid in $\cat{Met_*^p}$ is a
  $p$-subadditive commutative metric monoid.
		
  A morphism $f:(M,d,+_M,m_0) \to (N,\rho,+_N,n_0)$ of commutative
  monoids in $\cat{Met_*^p}$ is a morphism
  $f:(M,d,m_0) \to (N,\rho,n_0) \in \cat{Met_*^p}$ such that for all
  $a,b \in M$, $f(a+_M b) = f(a) +_N f(b)$ and $f(m_0) = n_0$.  That
  is, $f:(M,d,m_0) \to (N,\rho,n_0) \in \cat{Met_*}$ such that
  $f:(M,+_M,m_0) \to (N,+_N,n_0)$ is a monoid homomorphism.
		
  To see that $\cat{CMon(Met_*,\times_p,*)} = \cat{CMon(Met_*^p)}$,
  note that is suffices to require that $(M,d,m_0) \in \cat{Met_*}$,
  since the unit condition and $p$-subadditivity implies that for
  $a,b \in M$,
  $d(a,b) = d(a+m_0,m_0+b) \leq \norm{(d(a,m_0),d(m_0,b))}_p$.
\end{proof}

\subsection{Wasserstein distance}
\label{sec:wasserstein-distance}

In this section,
we introduce the $p$-Wasserstein distance $W_p$ on the space of diagrams $D(X,A)$ on a metric pair, and show that $(D(X,A),W_p)$, taken together with the monoid structure on $D(X,A)$, forms a $p$-subadditive commutative metric monoid.

Given a set pair $(X,A)$, recall that $D(X,A) = D(X)/D(A)$. As a special case, for a pointed set $(X,x_0)$, $D(X,x_0) = D(X)/D(x_0)$.

\begin{definition} \label{def:wasserstein}
  Let $p \in [1,\infty]$.
  Given a metric pair $(X,d,A)$ define $W_p[d,A]: D(X,A) \times D(X,A) \to [0,\infty]$ by
  \begin{equation*}
    W_p[d,A](x_1 + \cdots + x_m, x'_1 + \cdots + x'_n) = \inf \norm{ ( d(x_k,x'_{\sigma(k)}) )_{k=1}^{m+n} }_p,
  \end{equation*}
  where the infimum is taken over $x_{m+1}, \ldots, x_{m+n}, x'_{n+1}, \ldots x'_{n+m} \in A$ and $\sigma \in \Sigma_{m+n}$, where
  $\Sigma_{m+n}$ denotes the symmetric group on $m+n$ symbols.
\end{definition}

%
One may check that Definition~\ref{def:wasserstein} may be restated as follows.
\begin{lemma} \label{lem:Wp}
  \begin{equation*}
    W_p[d,A](x_1 + \cdots + x_m, x'_1 + \cdots + x'_n) = \min_{\sigma \in \Sigma_{m+n}} \norm{ ( d(x_k,x'_{\sigma(k)}) )_{k=1}^{m+n} }_p,
    \end{equation*}
  where $x_{m+1} = \cdots = x_{m+n} = A = x'_{n+1} = \cdots x'_{n+m}$ and $d(x,A) = \inf_{a \in A}d(x,a)$.
\end{lemma}
  For brevity, we will sometimes denote $W_p[d,A]$ by $W_p$ when this can lead to no confusion.
  
    	\begin{remark}\label{rem:quotient_iso}
  		Recall that for a metric pair $(X,d,A)$, $(X/A,\overline{d},A)$ denotes the pointed metric space obtained by collapsing $A$ to a point (see Definition \ref{def:quotient-pointed-metric-space}). Then $D(X,A)\cong D(X\setminus A) \cong D(X/A,A)$. Explicitly, we have monoid isomorphisms $\varphi:D(X\setminus A)\to D(X,A)$ and $\psi:D(X\setminus A)\to D(X/A,A)$ given by $x_1 + \dots + x_n \mapsto x_1 + \dots + x_n + D(A)$ and $x_1 + \dots + x_n \mapsto [x_1] + \dots + [x_n] + D(A)$, respectively. By Lemma~\ref{lem:natiso} below, we have that $(D(X,A),W_p[d,A], + ,0)$ and $(D(X/A,A),W_p[\overline{d}_p,A], + , 0)$ are isometrically isomorphic, and so we may pass between the settings of metric pairs and pointed metric spaces whenever convenient.
  	\end{remark}	

\begin{example} \label{ex:wasserstein}
  For the metric pair $(\R^2_{\leq},d,\Delta)$ or $(\overline{\R}^2_{\leq},d,\Delta)$ (Example~\ref{ex:half-plane}) and $p \in [1,\infty]$, $W_p[d,\Delta]$ is the $p$-Wasserstein distance on (finite) persistence diagrams.
  For the metric pair $(\Int(\R),d,\emptyset)$ and $d$ the length of the symmetric difference (Example~\ref{ex:intervals}), $W_1[d,\emptyset]$ is the barcode metric.
\end{example}

The following lemma verifies that $W_p$ is indeed a metric on $D(X,x_0)$.
\begin{lemma}
  Let $p \in [1,\infty]$.
  If $(X,d,x_0)$ is a pointed metric space then $(D(X,x_0),W_p,0)$ is a pointed metric space.  
\end{lemma}

\begin{proof}

    Point triviality and symmetry follow from the definition.
  To prove the triangle inequality, let
  $ \alpha = x_{1} + \dots + x_{n}, 
  \beta = x_{1}' + \dots + x_{m}',
  \gamma = x_{1}'' + \dots + x_{p}''$, be elements of $D(X,x_0)$. Let
  $r = n + m + p$
  and let $x_{n+1} = \cdots = x_r = x'_{m+1} = \cdots = x'_r = x''_{p+1} = \cdots = x''_r = x_0$.
  Let
  $\sigma$, $\tau\in S_{r}$ be permutations realizing
  $W_p(\alpha,\gamma)$, $W_p(\gamma,\beta)$, respectively. Let
  $\pi = \tau\circ\sigma\in S_r$. Then 

\begin{align*} W_p(\alpha,\beta)
    &\leq \elp{\big(d(x_k,x_{\pi(k)}')\big)_{k = 1}^r}{p}\\
    & \leq \elp{\big(d(x_k,x_{\sigma(k)}'')\big)_{k = 1}^r + \big(d(x_{\sigma(k)}'',x_{\pi(k)}')\big)_{k = 1}^r}{p}\\
    & \leq \elp{\big(d(x_k,x_{\sigma(k)}'')\big)_{k = 1}^r}{p} + \elp{\big(d(x_{\sigma(k)}'',x_{\pi(k)}')\big)_{k = 1}^r}{p}\\
    & = \elp{\big(d(x_k,x_{\sigma(k)}'')\big)_{k = 1}^r}{p} + \elp{\big(d(x_{\sigma(\sigma^{-1}(\ell))}'',x_{\pi(\sigma^{-1}(\ell))}')\big)_{\ell = 1}^r}{p} \\
    & = \elp{\big(d(x_k,x_{\sigma(k)}'')\big)_{k = 1}^r}{p} +
      \elp{\big(d(x_{\ell}'',x_{\tau(\ell)}')\big)_{\ell = 1}^r}{p} =
      W_p(\alpha,\gamma) + W_p(\gamma,\beta).\qedhere
  \end{align*}
\end{proof}

By Remark \ref{rem:quotient_iso}, the preceding lemma shows that $W_p[d,A]$ is a metric on $D(X,A)$ for any metric pair $(X,d,A)$ and any $p\in [1,\infty]$.

Next, we show that $W_p$ is $p$-subadditive.
\begin{lemma} \label{lem:p-subadditive}
  Let $p \in [1,\infty]$.
  If $(X,d,x_0)$ is a pointed metric space then $(D(X,x_0),W_p,+,0)$ is a $p$-subadditive commutative metric monoid and hence is an object in $\cat{CMon(Met_*^p)}$.
\end{lemma}

\begin{proof}

    Let $\alpha,\beta,\gamma,\delta \in D(X,x_0)$, where
  $\alpha = x_1+\cdots+x_m$, $\beta=x_{m+n+1}+\cdots+x_{m+n+p}$,
  $\gamma = x'_1+\cdots+x'_n$, $\delta=x'_{m+n+1}+\cdots+x_{m+n+q}$.
  We want to show that $W_p(\alpha+\beta,\gamma+\delta) \leq \norm{(W_p(\alpha,\gamma),W_p(\beta,\delta))}_p$.
  Let $x_{m+1}=\cdots=x_{m+n}=x_{m+n+p+1} = x_{m+n+p+q} =
  x'_{n+1}=\cdots=x'_{m+n}=x'_{m+n+q+1} = x'_{m+n+p+q} = x_0$.
Given $\sigma \in \Sigma_{m+n}$ and $\tau \in \Sigma_{p+q}$,  let $\sigma*\tau \in \Sigma_{m+n+p+q}$ be defined by $\sigma*\tau(i) := \sigma(i)$ if $i \leq m+n$ and $\sigma*\tau(i) := \tau(i)$ otherwise. Then

\begin{align*}
    \norm{(W_p(\alpha,\gamma),W_p(\beta,\delta))}_p &=
    \norm{ \left( \min_{\sigma \in \Sigma_{m+n}} 
    \norm{ ( d(x_i,x'_{\sigma(i)}) )_{i=1}^{m+n} }_p, 
    \min_{\tau \in \Sigma_{p+q}}
    \norm{ ( d(x_{m+n+i},x'_{m+n+\tau(i)}) ) }_p \right) }_p\\
    &= \min_{\sigma \in \Sigma_{m+n}} \min_{\tau \in \Sigma_{p+q}} 
      \norm{ (d(x_i,x'_{\sigma*\tau(i)}))_{i=1}^{m+n+p+q}}_p\\
    &\geq \min_{\pi \in \Sigma_{m+n+p+q}} 
      \norm{ (d(x_i,x'_{\pi(i)}))_{i=1}^{m+n+p+q}}_p
    = W_p(\alpha+\beta,\gamma+\delta).\qedhere
  \end{align*}
\end{proof}

\begin{lemma}\label{lem:i_is_isometry}
  Let $p \in [1,\infty]$.
  Let $(X,d,x_0) \in \cat{Met_*^p}$.
  Then the inclusion map $i:(X,d,x_0) \incl (D(X,x_0),W_p[d,x_0],0)$ is an isometry (and hence a metric map).
\end{lemma}

\begin{proof}
  Let $x,y \in X$.
  $W_p[d,x_0](x,y) = \min\left( d(x,y), \norm{(d(x,x_0),d(x_0,y))}_p \right) = d(x,y)$.
\end{proof}

The preceding lemma shows that if $(X,d,x_0)\in \cat{Met_*^p}$, then $i^*W_p = d$. On the other hand, if $(X,d,x_0)\in \cat{Met_*}$ but $d$ does not satisfy the $p$-strengthened inequality with respect to $x_0$, then the inclusion $i:(X,d,x_0)\hookrightarrow (D(X,x_0),W_p,x_0)$ is only guaranteed to be $1$-Lipschitz, but is not in general an isometry. This is one reason for working with $\cat{Met_*^p}$ as opposed to just $\cat{Met_*}$.

The following lemma shows that for a pointed metric map $f:(X,d_X,x_0)\to (Y,d_Y,y_0)$, the induced map $f_*:D(X,x_0)\to D(Y,y_0)$ is a metric map with respect to the Wasserstein distances.
\begin{lemma}
  Let $p \in [1,\infty]$.
  Given a pointed metric map $f:(X,d,x_0) \to (Y,d',y_0)$, the induced map $f_*:D(X,x_0) \to D(Y,y_0)$ is a morphism of $p$-subadditive commutative metric monoids $f_*:(D(X,x_0),W_p[d,x_0],+,0) \to (D(Y,y_0),W_p[d',y_0],+,0)$.
\end{lemma}

\begin{proof}

    By Definition~\ref{def:metric-monoid}, we need to show that $f_*:(D(X,x_0),W_p[d,x_0]) \to (D(Y,y_0),W_p[d',y_0])$ is a metric map and that
  $f_*:(D(X,x_0),+,0) \to (D(Y,y_0),+,0)$ is a monoid homomorphism.
  The latter is true by the definition of $f_*$ (Definition~\ref{def:persistence-diagram}).
  Let $x_1+\cdots+x_m, x'_1+\cdots+x'_n \in D(X,x_0)$.
  Let $x_{m+1} = \cdots = x_{m+n} = x_0 = x'_{n+1} = \cdots =x'_{n+m}$
  and thus $f(x_{m+1}) = \cdots = f(x_{m+n}) = y_0 = f(x'_{n+1}) = \cdots = f(x'_{n+m})$. Then

\begin{multline*}
    W_p[d',y_0](f(x_1)+\cdots+f(x_m),f(x'_1)+\cdots+f(x'_n))
  = \min_{\sigma \in \Sigma_{m+n}} \norm{ (d'(f(x_i),f(x'_{\sigma(i)})))_{i=1}^{m+n} }_p\\
  \leq \min_{\sigma \in \Sigma_{m+n}} \norm{ (d(x_i,x'_{\sigma(i)}))_{i=1}^{m+n} }_p
  = W_p[d,x_0](x_1+\cdots+x_m,x'_1+\cdots+x_n).\qedhere
\end{multline*}
\end{proof}

From the preceding lemmas, it is easy to see that the assignment that sends $(X,d,x_0)\in \cat{Met_*^p}$ to $(D(X,x_0),W_p,+,0)$ and that sends a pointed metric map $f:(X,d_X,x_0)\to (Y,d_Y,y_0)$ to the induced map $f_*:D(X,x_0)\to D(Y,y_0)$ is functorial.
\begin{definition} \label{def:Dp}
  Let $p \in [1,\infty]$.
  Let $D_p:\cat{Met_*^p} \to \cat{CMon(Met_*^p)}$ be the functor given by sending $(X,d,x_0)$ to $(D(X,x_0),W_p,+,0)$ and $f:(X,d,x_0) \to (Y,d',y_0)$ to 
  $f_*:(D(X,x_0),W_p[d,x_0],+,0) \to (D(Y,y_0),W_p[d',y_0],+,0)$.
\end{definition}

Recall that there is a forgetful functor $U_p:\cat{CMon(Met_*^p)}\to \cat{Met_*^p}$ given by sending $(M,d,+,0)\in \cat{CMon(Met_*^p)}$ to $(M,d,0)\in\cat{Met_*^p}$ (see Lemma \ref{lem:forgetful_functor_U} and Theorem~\ref{thm:CMon}).

\begin{theorem} \label{thm:Dp}
  Let $p \in [1,\infty]$.
  The forgetful functor $U_p:\cat{CMon(Met_*^p)} \to \cat{Met_*^p}$ has left adjoint $D_p$.
\end{theorem}

\begin{proof}
  Let $(X,d,x_0) \in \cat{Met_*^p}$, $(N,\rho,+,0) \in \cat{CMon(Met_*^p)}$ and $\varphi: (X,d,x_0) \to (N,\rho,0) \in \cat{Met_*^p}$.
  \begin{equation} \label{cd:3}
    \begin{tikzcd}
      (X,d,x_0) \ar[r,"i"] \ar[dr,swap,"\varphi"]& (D(X,x_0),W_p[d,x_0],0) \ar[d,dashed,"\tilde{\varphi}"]\\
      & (N,\rho,0)
    \end{tikzcd}
    \quad \quad \quad
    \begin{tikzcd}
      (D(X,x_0),W_p[d,x_0],+,0) \ar[d,dashed,"\exists ^^21 \tilde{\varphi}"] \\
      (N,\rho,+,0)
    \end{tikzcd}
  \end{equation}
  By the commutativity of the left hand side of \eqref{cd:3}, we have that for all $x \in X$, $\tilde{\varphi}(x) = \tilde{\varphi}(i(x)) = \varphi(x)$.
  For $\tilde{\varphi}$ to be a monoid homomorphism, we have
  $\tilde{\varphi}(x_1+\cdots+x_n) = \varphi(x_1) + \cdots + \varphi(x_n)$.
  Thus, if $\tilde{\varphi}$ exists it is unique. It remains to show that $\tilde{\varphi}$ is a metric map.

  Let $\alpha = x_1+\cdots+x_m \in D(X,x_0)$, $\beta = x'_1+\cdots+x'_n \in D(X,x_0)$ and let
  $x_{m+1}=\cdots=x_{m+n}= x_0 = x'_{n+1} = \cdots = x'_{n+m}$. Then
  \begin{multline*}
\rho(\tilde{\varphi}(\alpha),\tilde{\varphi}(\beta))
    = \rho(\varphi(x_1)+\cdots+\varphi(x_{m+n}),\varphi(x'_1)+\cdots+\varphi(x'_{m+n}))\\
    \leq \min_{\sigma \in \Sigma_{m+n}} \norm{ (\rho(\varphi(x_i),\varphi(x'_i)))_{i=1}^{m+n}}_p
    \leq \min_{\sigma \in \Sigma_{m+n}} \norm{ (d(x_i,x'_i))_{i=1}^{m+n}}_p
    = W_p[d,x_0](\alpha,\beta).
  \end{multline*}
Thus $((D(X,x_0),W_p[d,x_0],0),i)$ is a universal element, and the fact that $D_p$ is left adjoint to $U_p$ now follows from Lemma \ref{lem:adjoint}.
\end{proof}

The above constructions can also be formalized using metric pairs instead of pointed metric spaces, as we will now demonstrate.

\begin{definition} \label{def:U}
  Let $\overline{U}_p:\cat{CMon(Met_*^p)} \to \cat{Met_{pairs}}$ be the functor given by sending $(M,d,+,0)$ to $(M,d,\{0\})$ and $f:(M,d,+,0) \to (N,\rho,+,0)$ to $f:(M,d,\{0\}) \to (N,\rho,\{0\})$.
\end{definition}

Note that $\overline{U}_p$ is just the the composition $\cat{CMon(Met^p_*)}\xto{U_p}\cat{Met^p_*}\hookrightarrow\cat{Met_*}\hookrightarrow \cat{Met_{pairs}}$.

\begin{definition} \label{def:overlineDp}
  Let $\overline{D}_p:\cat{Met_{pairs}} \to \cat{CMon(Met_*^p)}$ be the functor given by sending $(X,d,A)$ to $(D(X,A),W_p[d,A],+,0)$ and $f:(X,d,A) \to (Y,d',B)$ to $f_*:(D(X,A),W_p[d,A],\allowbreak +,0) \to (D(Y,B),W_p[d',B],+,0)$.
\end{definition}

Recall the functors $Q:\cat{Met_{pairs}}\to \cat{Met_*}$ of Definition \ref{def:Q} and $S_p:\cat{Met_*}\to\cat{Met_*^p}$ of Definition \ref{def:Spq}.
We will show that $\overline{D}_p$ is the left adjoint of $\overline{U}_p$ (Theorem \ref{thm:overlineDp}). This will follow from the following lemma, which shows that $\overline{D}_p$ is naturally isomorphic to the composition $D_pS_pQ$, together with the fact that each of $D_p,S_p,Q$ has a right adjoint, and the composition of these right adjoints is precisely $\overline{U}_p$.

\begin{lemma}\label{lem:natiso}
	The functors $\overline{D}_p:\cat{Met_{pairs}}\to \cat{CMon(Met_*^p)}$ and $D_pS_pQ:\cat{Met_{pairs}}\to \cat{CMon(Met_*^p)}$ are naturally isomorphic.
\end{lemma}

\begin{proof}

    Let $(X,d,A) \in \cat{Met_{pairs}}$.
    Then $\overline{D}_p(X,d,A) = (D(X,A),W_p[d,A],+,0)$ and
    $D_pS_pQ(X,d,A) = (D(X/A,A),W_p[\overline{d}_p,A],+,0)$.
    Recall (see Remark~\ref{rem:quotient_iso}) that we have monoid isomorphisms $\varphi = \varphi_X:D(X \setminus A) \to D(X,A)$ and $\psi = \psi_X:D(X \setminus A) \to D(X/A,A)$. Let $\eta = \eta_X:D(X,A) \to D(X/A,A)$ be the composite monoid isomorphism $\psi_X \varphi_X^{-1}$. Explicitly, $\eta_X$ is given by $x_1+\dots + x_n + D(A)\mapsto [x_1] + \dots + [x_n] + D(A)$.
    We will show that $\eta_X$ is a isometry.
    Let $x_1+\cdots+x_m,x'_1+\cdots+x'_n \in D(X \setminus A)$.
    Denote these elements by $\alpha$ and $\alpha'$, respectively.
    Let $x_{m+1}= \cdots = x_{m+n} = x'_{n+1} = x'_{n+m} = A$.
    By Lemma~\ref{lem:Wp},

  \begin{equation*}
      W_p[d,A](\varphi\alpha,\varphi\alpha') = \min_{\sigma \in \Sigma_{m+n}}\norm{(d(x_i,x_{\sigma(i)}'))_{i=1}^{m+n}}_p.
    \end{equation*}
    On the other hand, by Lemma~\ref{lem:overlineDp},
    \begin{align*}
      W_p[\overline{d}_p,A](\psi\alpha,\psi\alpha')
    &= \min_{\sigma \in \Sigma_{m+n}}\norm{(\overline{d}_p(x_i,x_{\sigma(i)}'))_{i=1}^{m+n}}_p\\
      &= \min_{\sigma \in \Sigma_{m+n}}\norm{(\min(d(x_i,x_{\sigma(i)}'),\norm{d(x_i,A),d(A,x_{\sigma(i)})}_p))_{i=1}^{m+n}}_p.
    \end{align*}

      Therefore
      $W_p[\overline{d}_p,A](\psi\alpha,\psi\alpha')\leq  W_p[d,A](\varphi\alpha,\varphi\alpha')$.
      On the other hand, let $\sigma_0\in \Sigma_{n+m}$ be a permutation such that $W_p[\overline{d}_p,A](\psi\alpha,\psi\alpha') =
      \norm{(\min(d(x_i,x_{\sigma_0(i)}'),\norm{d(x_i,A),d(A,x_{\sigma_0(i)})}_p))_{i=1}^{m+n}}_p$.
If there is an index $i$ with $\norm{(d(x_i,A),d(A,x'_{\sigma_0(i)}))}_p < d(x_i,x_{\sigma_0(i)}')$ then we may choose a new permutation $\sigma_1$ with $x'_{\sigma_1(i)} = A$ and $x_{\sigma_1^{-1}\sigma_0(i)} = A$ and which is otherwise the same as $\sigma_0$.
By induction we remove all such indices to obtain
a new permutation $\sigma_0'\in \Sigma_{n+m}$ with 
      \[\lVert(d(x_i,x_{\sigma_0'(i)}'))_{i=1}^{m+n}\rVert_p = \norm{(\overline{d}_p(x_i,x_{\sigma_0(i)}'))_{i=1}^{m+n}}_p.\]
      Therefore
$W_p[\overline{d}_p,A](\psi\alpha,\psi\alpha') = W_p[d,A](\varphi\alpha,\varphi\alpha')$. Then 

	\begin{multline*}
			W_p[\overline{d}_p,A](\eta\alpha,\eta\alpha') = W_p[\overline{d}_p,A](\psi\varphi^{-1}\alpha,\psi\varphi^{-1}\alpha')\\
			= W_p[d,A](\varphi\varphi^{-1}\alpha,\varphi\varphi^{-1}\alpha') = W_p[d,A](\alpha,\alpha').
		\end{multline*}
    Thus we have an isomorphism $\eta_X: \overline{D}_p(X,d,A) \to D_pS_pQ(X,d,A)$.
 
  To see that these isomorphisms are natural, let $f:(X,d,A)\to (Y,d',B)$ be a morphism in $\cat{Met_{pairs}}$. The map $\overline{D}_pf = f_*:D(X,A)\to D(Y,B)$ is given by $x_1 + \dots + x_n + D(A)\mapsto f(x_1) + \dots + f(x_n) + D(B)$, while the map $D_pS_pQf:D(X/A,A)\to D(Y/B,B)$ is given by $[x_1] + \dots + [x_n] + D(A)\mapsto [f(x_1)] + \dots + [f(x_n)] + D(B)$. Thus $(D_pS_pQf)\eta_X = \eta_Y\overline{D}_pf$ and hence the maps $\eta_X$ assemble into a natural isomorphism $\eta:\overline{D}_p\Rightarrow D_pS_pQ$.
\end{proof}

\begin{theorem}\label{thm:overlineDp}
  The forgetful functor $\overline{U}_p:\cat{CMon(Met_*^p)} \to \cat{Met_{pairs}}$ has left adjoint $\overline{D}_p$.
\end{theorem}

\begin{proof} 
  Since $\overline{U}_p$ is given by the composition $\cat{CMon(Met_*^p)} \xto{U_p} \cat{Met_*^p} \incl \cat{Met_*} \incl \cat{Met_{pairs}}$,
  by Theorems \ref{thm:Q}, \ref{thm:Spq}, and \ref{thm:Dp}, it has left adjoint the composite $D_p S_p Q$.
  By Lemma \ref{lem:natiso}, $\overline{D}_p\cong D_pS_pQ$ and hence $\overline{D}_p$ is, up to natural isomorphism, the left adjoint of $\overline{U}_p$.
\end{proof}

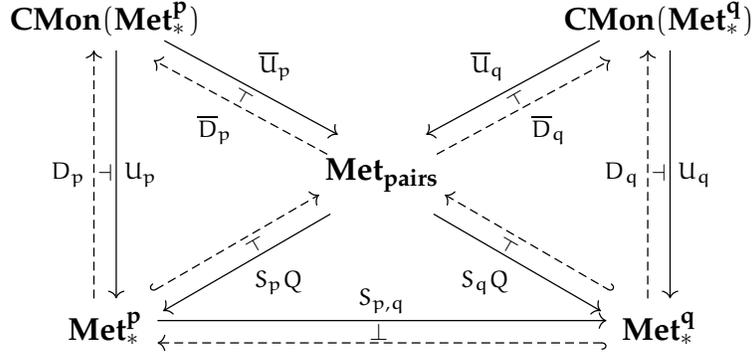
\begin{figure}
	\[\begin{tikzcd}[column sep=3em,row sep=3em]
		\cat{CMon(Met_*^p)} \ar["U_p"{name=Up,right},shift left = .75ex]{dd}
		\ar[""{name=Upbar,right},shift left=.75ex]{dr}[xshift = 0ex]{\overline{U}_p}
		& & \cat{CMon(Met_*^q)}
		\ar[""{name=Uq,right},shift left = .75ex]{dd}{U_q}
		\ar[""{name=Uqbar,right},shift right = .75ex,swap]{dl}{\overline{U}_q}\\
		& \cat{Met_{pairs}}
		 \ar[""{name=Iq,left},hookleftarrow, dashed,shift left = .75ex]{dr}
		 \ar[""{name=Dpbar,left},dashed,shift left = .75ex]{ul}[xshift = 0ex]{\overline{D}_p}
		 \ar[""{name=Dqbar,left},dashed,shift right = .75ex, swap]{ur}{\overline{D}_q}
		 \ar[""{name=SpQ,left},shift left = .75ex]{dl}{S_pQ}
		 \ar[""{name=SqQ,right},swap,shift right = .75ex]{dr}{S_qQ}&\\
		\cat{Met_*^p}  
		\ar[""{name=Ip,right},hookrightarrow,dashed,shift left = .75ex]{ur}{}
		\ar[""{name=Spq,left},shift left = .75ex]{rr}{S_{p,q}}
		\ar["D_p"{name=Dp,left},dashed,shift left = .75ex]{uu}
		& & \cat{Met_*^q}
		\ar[""{name=Ipq,left}, hookrightarrow,dashed,shift left = .75ex]{ll}
		\ar[""{name=Dq,left}, dashed, shift left = .75ex]{uu}{D_q}
		\ar[from=Up,to=Dp,draw=none,shift right = 1.25ex]{}{\dashv}
		\ar[from=Upbar,to=Dpbar,draw=none]{}[rotate=60,xshift=-1.1ex,yshift=1.25ex]{\dashv}
		\ar[from=Uqbar,to=Dqbar,draw=none]{}[rotate=120,xshift=-1.1ex,yshift=1.25ex]{\dashv}\ar[from=Uq,to=Dq,draw=none,shift right = 1.25ex]{}{\dashv}
		\ar[from=Ip,to=SpQ,draw=none]{}[rotate=120,xshift=-1.1ex,yshift=0ex]{\dashv}	\ar[from=Iq,to=SqQ,draw=none]{}[rotate=60,xshift=-1.1ex,yshift=0ex]{\dashv}
		\ar[from=Ipq,to=Spq,draw=none]{}[xshift=1.5ex]{\perp}
	\end{tikzcd}
	\]
	\caption{For $1\leq p\leq q\leq \infty$, the relationships between the categories $\cat{Met_*^p}$, $\cat{Met_*^q}$, $\cat{CMon(Met_*^p)}$, $\cat{CMon(Met_*^q)}$, and $\cat{Met_{pairs}}$ via the functors $U_r, D_r, \overline{U}_r, \overline{D}_r$, $S_r$ ($r \in \{p,q\}$), and $Q, S_{p,q}$.
            The solid triangles commute and the dashed triangles commute.}
        \label{fig:big_comm_diag}
\end{figure}

The relationship between the forgetful functors $U_p,\overline{U}_p$, the free functors $D_p,\overline{D}_p$, the quotient functors $Q$, and the functors $S_p, S_{p,q}$ is summarized in Figure \ref{fig:big_comm_diag}.

\section{Applications}
\label{sec:applications} In this section, we give several applications of universality. The first application shows that for a pointed metric space $(X,d,x_0)$, $W_p$ is the largest $p$-subadditive metric on $D(X,x_0)$ which in some sense extends the metric $d$. This result implies an abstract form of converse stability from which we derive converse stability-type results in various settings. As a second application, we show how universality can be used to derive the correct form of Kantorovich-Rubinstein duality for persistence diagrams.

\subsection{Maximality of the Wasserstein distances}
\label{sec:maximality}
The following theorem shows that $W_p$ is the largest $p$-subadditive metric extending the underlying metric.
\begin{theorem} \label{thm:maximal}
	Let $p\in [1,\infty]$ and let $(X,d,A)$ be a metric pair. Then $W_p[d,A]$ is the largest $p$-subadditive metric $\rho$ on $D(X,A)$ satisfying $i^*\rho = d_p$.
\end{theorem}

\begin{proof}
	Suppose that $\rho$ is a $p$-subadditive metric on $D(X,A)$ with $i^*\rho = d_p$. Then $(D(X,A),\rho,+,0)\in \cat{CMon(Met_*^p)}$ and $i:(X,d,A)\hookrightarrow (D(X,A), \rho, +, 0)$ is $1$-Lipschitz. By Theorem \ref{thm:overlineDp}, there is a unique $1$-Lipschitz map $\tilde{i}:(D(X,A),W_p[d,A],0)\to (D(X,A),\rho,0)$, and hence $\rho\leq W_p$.
\end{proof}

If $d$ satisfies the $p$-strengthened triangle inequality then $d = d_p$, and so we immediately obtain the following.

\begin{corollary}\label{cor:maximal}
	Let $p\in [1,\infty]$ and let $(X,d,x_0)\in \cat{Met_*^p}$. Then $W_p[d,x_0]$ is the largest $p$-subadditive metric $\rho$ on $D(X,x_0)$ satisfying $i^*\rho = d$.
\end{corollary}

As another application of universality, we show that for a pointed metric space $(X,d,x_0)$, $W_p[d,x_0] = W_p[d_p,x_0]$.
\begin{corollary}\label{cor:Wp[dp]equalWp[d]}
	Let $(X,d,x_0)\in \cat{Met_*}$. Then $W_p[d,x_0] = W_p[d_p,x_0]$. 
	\begin{proof}
		Note that $i^*W_p[d,x_0] = i^*W_p[d_p,x_0] = d_p$. Since $W_p[d,x_0]$ and $W_p[d_p,x_0]$ are both $p$-subadditive, we have $W_p[d,x_0] = W_p[d_p,x_0]$ by Theorem \ref{thm:maximal}.
	\end{proof}
\end{corollary}

The preceding corollary justifies our use of the categories $\cat{Met_*^p}$. By Lemma \ref{lem:i_is_isometry}, if $(X,d,x_0)\in \cat{Met_*^p}$, then $(X,d)$ embeds into $(D(X,x_0),W_p[d,x_0],0)$. This corollary shows that we can always pass to $d_p$ without changing the Wasserstein distance, and so it suffices to work in $\cat{Met_*^p}$.

\subsection{Converse stability}
\label{sec:converse}

We show that certain converse stability theorems follow from our results. The following is a completely formal (or ``soft''~\cite{bubenik2015metrics}) converse stability result, from which specific converse stability theorems follow.

\begin{theorem}[Abstract Converse Stability]\label{thm:stability}
	Fix $p\in [1,\infty]$. Let $(X,x_0)$ be a pointed set and let $\rho$ be a $p$-subadditive metric on $D(X,x_0)$. Then $\rho\leq W_p[i^*\rho,x_0]$. 
\end{theorem}

\begin{proof}
	Let $d = i^*\rho$ and consider the pointed metric space $(X,d,x_0)$. Since $\rho$ is $p$-subadditive by assumption and $i^*\rho = d$ by definition, the result immediately follows from Corollary \ref{cor:maximal}.
\end{proof}

\begin{example}[Converse Algebraic Stability]
  \label{ex:converse_algebraic_stability}
  (See \cite{lesnick2015theory} for a version of this result that applies to all pointwise finite dimensional persistence modules.)
%
  Let $\mathbf{Vect}(K)^{\R}_{\textup{Fin}}$ denote the monoid of isomorphism classes of persistence modules which decompose as a finite direct sum of interval modules.
  We can identify $\mathbf{Vect}(K)^{\R}_{\textup{Fin}}$ with
  $D(\Int(\R),\emptyset)$
  via the map
that sends a direct sum of interval modules to the corresponding formal sum of intervals.
Equip $\mathbf{Vect}(K)^{\R}_{\textup{Fin}}$ with the interleaving distance $d_I$~\cite{chazal_cohen-steiner_glisse_guibas_oudot_2009,lesnick2015theory,bubenik2014categorification}. Note that $d_I$ is $\infty$-subadditive. Indeed, if $(\varphi,\psi)$ is an $\epsilon$-interleaving between $M$ and $N$ and $(\varphi',\psi')$ is an $\eta$-interleaving between $M'$ and $N'$, then $(\varphi \oplus \varphi',\psi \oplus \psi')$ is a $\max(\epsilon,\eta)$-interleaving between $M\oplus M'$ and $N\oplus N'$.
Note that the interleaving distance for interval modules
is $(d_H)_\infty$, the $\infty$-strengthening of the Hausdorff distance with respect to $\emptyset$. That is, $i^*d_I = (d_H)_\infty$, and hence by Theorem~\ref{thm:stability} and Corollary \ref{cor:Wp[dp]equalWp[d]},
$d_I \leq  W_\infty[(d_H)_\infty,\emptyset] = W_\infty[d_H,\emptyset]$.

For a second version of this result, introduce an equivalence relation on
$\mathbf{Vect}(K)^{\R}_{\textup{Fin}}$
given by $M \sim N$ if $d_I(M,N)=0$.
Then we can identify 
$\mathbf{Vect}(K)^{\R}_{\textup{Fin}}/\!\sim$ with $D(\overline{\R}^2_{\leq}, \Delta)$
via the map that sends a direct sum of interval modules to the corresponding persistence diagram.
By Theorem~\ref{thm:stability},
$d_I \leq W_{\infty}[d,\Delta]$, where $d$ is the $\ell^{\infty}$-distance.
The distance $W_{\infty}[d,\Delta]$ is the \emph{bottleneck distance}.
Note that this distance restricted to the images of interval modules is $d_{\infty}$, the $\infty$-strengthening of the $\ell^{\infty}$-distance.
\end{example}

  \begin{example}[Converse Algebraic Stability Theorem for generalized persistence modules]
    Consider generalized persistence modules
    $M: \mathbf{P} \to \mathbf{A}$.
    If $\mathbf{P}$ is equipped with certain additional structure, such as a \emph{subadditive projection on translations} or a \emph{superadditive family of translations}, then $\mathbf{A}^\mathbf{P}$ can be equipped with an interleaving distance $d_I$~\cite{bubenik2015metrics}.
  As in Example~\ref{ex:converse_algebraic_stability}, the interleaving distance $d_I$ is $\infty$-subadditive.
Let $\Ind$ be a set of indecomposable generalized persistence modules in 
$\mathbf{A}^{\mathbf{P}}$ with basepoint the zero module $0$.
Then we have the set of generalized barcodes $D(\Ind,0)$. 
There is a bijection from the set of isomorphism classes of generalized persistence modules in $\mathbf{A}^{\mathbf{P}}$ that are isomorphic to a finite direct sum of elements of $\Ind$ to $D(\Ind,0)$, which sends direct sums to formal sums.
By Theorem~\ref{thm:stability},
$d_I \leq W_{\infty}[d_I,0]$, where the latter is also called the bottleneck distance~\cite{bjerkevik2016stability}.
%
\end{example}

\subsection{Kantorovich-Rubinstein duality}
\label{sec:duality}

The classical Kantorovich-Rubinstein duality theorem says that the classical $1$-Wasserstein distance $w_1(\mu,\nu)$ between probability measures $\mu$ and $\nu$ on a complete and separable metric space $(X,d)$ is equal to $\sup\int_{X}f d(\mu-\nu)$, where the supremum is taken over all $1$-Lipschitz functions. A version of Kantorovich-Rubinstein duality holds for persistence diagrams as well. We will show that \[\textstyle W_1(\sum_{i = 1}^n a_i, \sum_{j = 1}^m b_j) = \sup\{\sum_{i = 1}^nk(a_i) - \sum_{j = 1}^mk(b_j) \},\] where now the supremum is taken over all $1$-Lipschitz functions $k:X\to \R$ with $k(x_0) = 0$.
	
To motivate the form that Kantorovich-Rubinstein duality takes for persistence diagrams, we first show how Theorem \ref{thm:2} can be used to derive the inequality 
\[\textstyle W_1(\sum_{i = 1}^n a_i, \sum_{j = 1}^m b_j) \geq \sup\{\sum_{i = 1}^nk(a_i) - \sum_{j = 1}^m k(b_j)\}.\]
 We will then use the classical Kantorovich-Rubinstein duality theorem to show that this is in fact an equality.

Let $(X,d,x_0)\in \cat{Met_*}$ and consider the commutative metric monoid $(\R,|\cdot|,+,0)$, where $|\cdot|$ denotes the metric induced by absolute value. The inequality $|(a+b)-(c+d)| \leq |a-c| + |b-d|$ implies that $(\R,|\cdot|,+,0) \in \cat{CMon(Met_*)}$.
Let $h:X\to \R$ be a $1$-Lipschitz map.
Define $k:X\to \R$ by $k(x) = h(x) - h(x_0)$ for all $x\in X$.
Then $k(x_0) = 0$ and $|k(x)-k(y)| = |h(x)-h(x_0)- (h(y)- h(x_0))| = |h(x)-h(y)| \leq d(x,y)$
so that $k$ is a pointed metric map.
By Theorem \ref{thm:2}, there is a unique morphism of $1$-subadditive commutative metric monoids $\tilde{k}:(D(X,x_0),\allowbreak W_1,+,0)\to (\R,|\cdot|,+,0)$
such that $\tilde{k}\circ i = k$.
Explicitly, $\tilde{k}$ is given by $\sum_i c_i \mapsto \sum_i k(c_i)$ for $c_i\in X$. Then for $\alpha = \sum_{i = 1}^n a_i$, $\beta = \sum_{j = 1}^m b_j \in D(X,x_0)$ with $m \geq n$, we have
$|\tilde{k}(\alpha) - \tilde{k}(\beta)| = \left|\tilde{k}\left(\sum_{i = 1}^n a_i\right) - \tilde{k}\left(\sum_{j = 1}^mb_j\right)\right|
= \left|\sum_{i = 1}^nk(a_i) - \sum_{j = 1}^m k(b_j)\right| = \left|\sum_{i = 1}^nh(a_i) - \sum_{j = 1}^mh(b_j) +(m-n)h(x_0)\right|$.
Since $\tilde{k}$ is $1$-Lipschitz, we obtain the inequality
$\left|\sum_{i = 1}^nh(a_i) - \sum_{j = 1}^mh(b_j) +(m-n)h(x_0)\right| \leq W_1\big(\sum_{i = 1}^n a_i, \sum_{j = 1}^m b_j\big)$.
Therefore
\begin{equation}\label{eq:KR_duality_inequality}
  \small
\sup\bigg\{\sum_{i = 1}^nh(a_i) - \sum_{j = 1}^mh(b_j) + (m-n)h(x_0) \ \bigg| \ h:X\to \R,\ 1\text{-Lipschitz}\bigg\} \leq W_1\bigg(\sum_{i = 1}^n a_i, \sum_{j = 1}^m b_j\bigg),
\end{equation}
or equivalently,
\begin{equation*}
  \small
\sup\bigg\{\sum_{i = 1}^nk(a_i) - \sum_{j = 1}^mk(b_j) \ \bigg| \ k:X\to \R, \, k(x_0) = 0, \, \ 1\text{-Lipschitz}\bigg\} \leq W_1\bigg(\sum_{i = 1}^n a_i, \sum_{j = 1}^m b_j\bigg).
\end{equation*}

%
To see that this inequality is in fact an equality, consider $\alpha = a_1 + \dots + a_n,\beta = b_1 + \dots + b_m\in D(X,x_0)$ and consider the classical $1$-Wasserstein distance $w_1(\sum_{i = 1}^n \delta_{a_i} + (r-n)\delta_{x_0}, \sum_{i = 1}^m \delta_{b_j} + (r-m)\delta_{x_0})$, where $\delta_x$ is the Dirac measure at $x$ and $r=m+n$. It is known that for 
sums of Dirac measures 
the computation of $w_1$ is equivalent to the \emph{linear assignment problem}. In other words, letting $\tilde{\alpha} = \sum_{i = 1}^n \delta_{a_i} + (r-n)\delta_{x_0}$ and $\tilde{\beta} = \sum_{i = 1}^m \delta_{b_j} + (r-m)\delta_{x_0}$, we have
	\[w_1\left(\tilde{\alpha}, \tilde{\beta}\right) = \min_{\sigma\in S_{n+m}}\elp{(d(a_i,b_{\sigma(i)})_{i = 1}^{n+m})}{1} = W_1[d,x_0](\alpha,\beta),\]
	where $a_{n+1} = \dots = {a_{n+m}} = b_{m+1} = \dots = b_{n+m} = x_0$.
	It follows from Kantorovich-Rubinstein duality for measures \cite{villani2003topics} that 
	\begin{equation}\label{eq:KR_duality_diracs} 
		\sup\left\{\int_Xhd(\tilde{\alpha}-\tilde{\beta}) \ | \ h:X\to \R, \ 1\text{-Lipschitz}\right\} = w_1(\tilde{\alpha},\tilde{\beta}) = W_1[d,x_0](\alpha,\beta).
	\end{equation} 
Since the left-hand side of \eqref{eq:KR_duality_diracs} is precisely the left-hand side of \eqref{eq:KR_duality_inequality}, we obtain the desired equality.

For a direct proof of Kantorovich-Rubinstein duality in this
setting using linear programming see \cite[Appendix
C]{be:universality}.


\subsection*{Acknowledgments}

This research was partially supported by the Southeast Center for Mathematics and Biology, an NSF-Simons Research Center for Mathematics of Complex Biological Systems, under National Science Foundation Grant No. DMS-1764406 and Simons Foundation Grant No. 594594.
This material is based upon work supported by, or in part by, the Army Research Laboratory and the Army Research Office under contract/grant number
W911NF-18-1-0307.

%

\newcommand{\etalchar}[1]{$^{#1}$}

\end{document}